\documentclass[a4paper,12pt,twoside]{article}

\usepackage[body={177mm,240mm}]{geometry}
\usepackage{a1}
\usepackage[english]{babel}
\usepackage[latin1]{inputenc} 
\usepackage{amsfonts,amssymb}
\usepackage{amsmath}
\usepackage{graphicx}
\usepackage{float}

\usepackage{makeidx}
\makeindex

\def\mapright#1#2#3{\smash{\mathop{\hbox to
#3{\rightarrowfill}}\limits^{#1}_{#2}}}

\def\mapleft#1#2#3{\smash{\mathop{\hbox to
#3{\leftarrowfill}}\limits^{#1}_{#2}}}

\def\mapright#1#2{\smash{\mathop{\hbox to 0.90cm{\rightarrowfill}}\limits^{#1}_{#2}}}
\def\mapleft#1#2{\smash{\mathop{\hbox to 0.90cm{\leftarrowfill}}\limits^{#1}_{#2}}}

\def\mapleftright#1#2{\smash{\mathop{\hbox to 0.80cm{\leftarrowfill \rightarrowfill}}\limits^{#1}_{#2}}}

\title{All the shapes of spaces: a census of small 3-manifolds
} 
\author{Sóstenes L. Lins and Lauro D. Lins}

\date{\today}

\begin{document}

\maketitle
\begin{abstract}

In this work we present a complete (no misses, no duplicates) census for closed, 
connected, orientable and
prime 3-manifolds induced by plane graphs with a 
bipartition of its edge set (blinks) up to $k=9$ edges. 
Blinks form a universal encoding for such manifolds. In fact, each such a manifold
is a subtle class of blinks, \cite{lins2013B}. Blinks are in 1-1 correpondence with
{\em blackboard framed links}, \cite {kauffman1991knots, kauffman1994tlr}
We hope that this census becomes as useful for the study of concrete examples of 3-manifolds
as the tables of knots are in the study of knots and links.

\end{abstract}

\section{Introduction} 

Similar census of 
3-manifolds homeomorphism classes encoded by blinks up to $k=10, 11, 12, \ldots$ edges
are, in principle, possible from the theory developed here.
The limit is the technological state of the art. The algorithms are fully parallelizable
and this limit seems far away. This work brings to light a new technique for obtaining
drawings of links with curls and blinks with loops and pendant edges. These decorative
objects in the usual theory of links become fundamental in this work. Links are drawn so
as to deterministically minimize the number of right angles turns in a grid. 
This comes as an application of
network flow optimization theory. Another contribution of this work is in finding good examples
to test the splendid general theory of (we are focused at closed oriented) 3-manifold
which is flourishing in this era post-Perelman.
We have done this by recently posting in the arXiv two papers with 2+11=13 
classes of 3-manifolds left open
in L. Lins' Thesis under the supervision of S. Lins, \cite{lins2007blink}. 
The first challenge, containing 2 classes was solved quickly by various 
researchers in the 3-manifold community. 
The second, with tougher challenges found by BLINK, having 11 classes of
3-manifolds is about of being solved by a combination of the softwares:
SnapPy/GAP/Sage. In the final section
of this paper we explain why these challenges should be taken seriously.

\subsection{Motivation}
After presenting some instances of closed 3-manifolds,
P. Alexandroff says in the English translation (1961) of 
his joint work with D. Hilbert \cite{alexandrov1961elementary}, first
published (1932) in German, \cite{alexandroff1932einfachste}:
{\em ``These few examples will suffice. Let it be remarked here that, at present,
in contrast with the two-dimensional case, the problem of enumerating the 
topological types of manifolds of three and more dimensions is in an apparently 
hopeless state. We are not only far removed from the solution, but even from the first step
toward a solution, a plausible conjecture''.}

John Hempel in his book (1976) {\em 3-Manifolds} \cite{hempel1976} writes 
at the opening of Section 15, entitled {Open Problems:}
{\em ``The ultimate goal of the theory would be in providing solutions to:
 {\em The homeomorphism problem}: provide an {\em effective} procedure
 for determining whether two {\em given} 3-manifolds are homeomorphic, together with
 {\em The classification problem}: {\em effectively} generate a list containing
 exactly one 3-manifold from each homeomorphism class.''}

 Along the years we were continuously motivated by the above quotations.

\subsection{A surprising new 1-1 correspondence:
$\frac{3-manifolds}{homeomorphisms}$ = $\frac{blinks}{coin \ calculus}$}
\label{subsec:surprise}
\begin{figure}[H]
\begin{center}
\includegraphics[width=15cm] {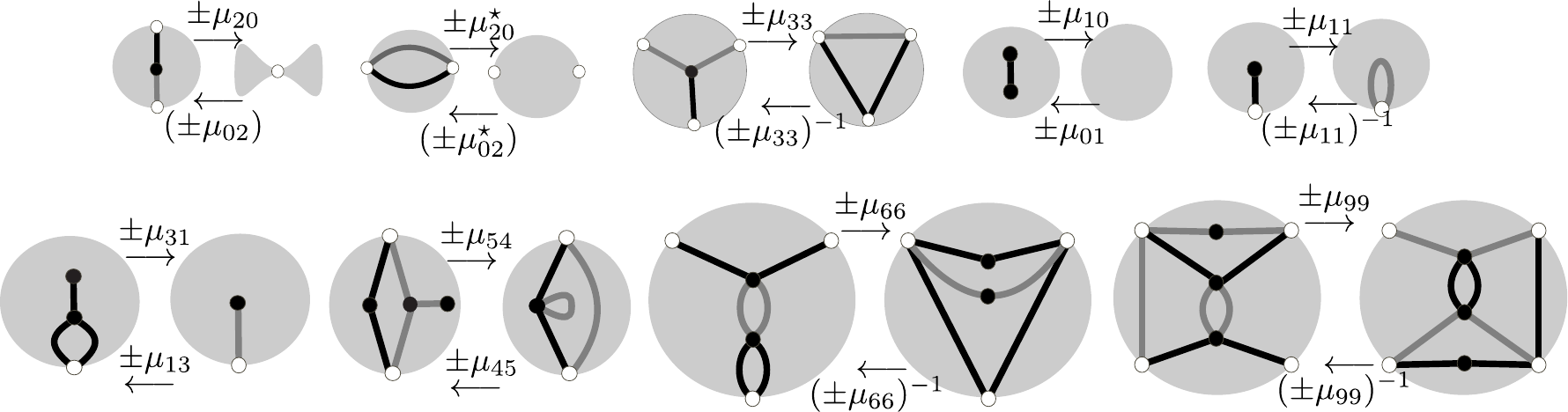}
\caption{\sf A 36-move, 36-coins (not all distinct) reduced (sufficient) coin calculus:
a coin is a sub-blink which lies in a disk (or in a pinched disk, in the case of the right coin
at the right of the moves $\pm \mu_{20}$). The complementary sub-blink is completely
arbitrary, provided the intersection with the coins are vertices contained in the
set of attachment vertices of the coins, each represented by a hollow small circle in the
boundary of the coin, as shown in Fig. \ref{fig:reducedblinkcalculus}.}
\label{fig:reducedblinkcalculus}
\end{center}
\end{figure}

A {\em blink} is a finite plane graph with an edge bipartition.
Any closed, connected and oriented 3-manifold is induced by some blink. Even though this object
has been around since 1994 when it was introduced in the joint research 
monography of L. Kauffman and
S. Lins, \cite{kauffman1994tlr}, the fact that they encode oriented closed 
3-manifolds remains basically unkown.
This is about to change because as a consequence, \cite{lins2013B}, of a recent
result of B. Martelli, \cite {martelli2011finite,martelli2012finite}, each such a 
3-manifold becomes a subtle equivalence class of blinks.
In \cite{lins2013B}, a new calculus, this time on blinks, 
named {\em the coin calculus} with 
8 types of local moves each applied to 8 types of related sub-blinks (named coins) are shown 
to capture the essence of homeomorphism between 3-manifolds, in the sense of 
Theorem \ref{theo:incredible}. A sufficent reduced coin calculus with 36 moves and 36 blinks 
is shown in Fig. \ref{fig:reducedblinkcalculus}.

\noindent
Throughout this work {\em 3-manifold} means a closed, oriented and connected 3-manifold.
\begin{theorem}
\label{theo:incredible} 
The class of closed, oriented and connected 3-manifolds modulo homeomorphism is equal to 
the class of blinks modulo the coin calculus of Fig. \ref{fig:reducedblinkcalculus}.
\end{theorem}
The meaning of the theorem is that there exists a bijection 
$$\beta: \ \ \frac{{3\hspace{-1mm}-\hspace{-1mm}manifolds}}{{homeomorphisms}}\ \  
\longrightarrow\ \ \  \ \frac{blinks}{coin \ calculus}.$$
If the manifold is given by a framed link, $\beta$ is obtained by a linear algorithm.
However, if it is given by a gem, by a triangulation, by a special spine 
\cite{matveev2007algorithmic},  
or by a Heegaard diagram, then a
polynomial algorithm to find $\beta$ seems to be an open problem. However, 
rescuing the situation there exists a recent work of S. Lins and his former student
R. Machado, which is reported in a sequence of 3-papers posted in the arXiv 
proving that there exists an 
$O(n^2)$-algorithm to go from a {\em resoluble} gem 
to a blink inducing the same 3-manifold. The definition of {\em resoluble gems} and
the status of this theory is the subject of Subsection \ref{subsec:on2algo}.
We antecipated that it solves the general problem up to a conjecture which is 
most certainly true. In the realm of the blinks and gems 
in the data basis of BLINK we have an $O(n^2)$-algorithm
to obtain a framed link inducing the same 3-manifold as an arbitrarily given gem.

\begin{figure}[H]
\begin{center}
\includegraphics[width=11.3cm] {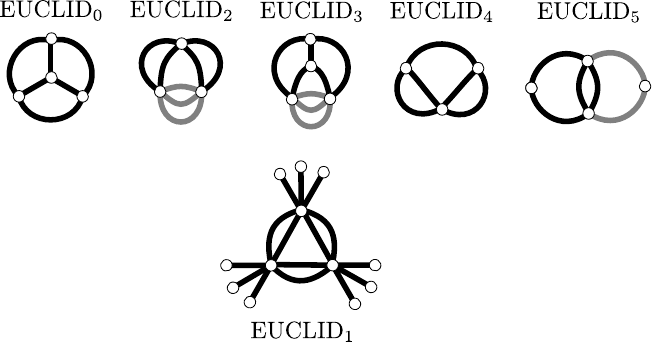}
\caption{\sf The discovery of a blink for EUCLID$_1$, solving an
open problem left in page 117 of \cite{lins2007blink}:
of the six euclidean 3-manifolds only EUCLID$_1$ did not have a blink presentation. 
It was necessary to develop and understanding some deep geometric properties of gems, 
to find this blink
by the theory developed in \cite{linsmachadoA2012,linsmachadoB2012,linsmachadoC2012}.
The new blink has more than twice the number of edges of the other blinks in the euclidean family.
No wonder it was so difficult to find. The new blink is correct: R. Machado applied to its
associated blackboard framed link the linear algorithm (given in Fig. \ref{fig:12to8move})
to produce the canonical gem $G$ inducing the same 3-manifold. The simplifying procedure of BLINK
was used to go from $G$ to the {\em superattractor}, (the
unique minimum gem with 24 vertices 
inducing $|\rm{EUCLID}_1|$), (\cite{lins1995gca}) of EUCLID$_1$.
}
\label{fig:euclideans}
\end{center}
\end{figure}  
 
\subsection{A complete duplicate free census of 9-small 3-manifolds}
\begin{figure}[H]
\begin{center}
\includegraphics[width=16.3cm] {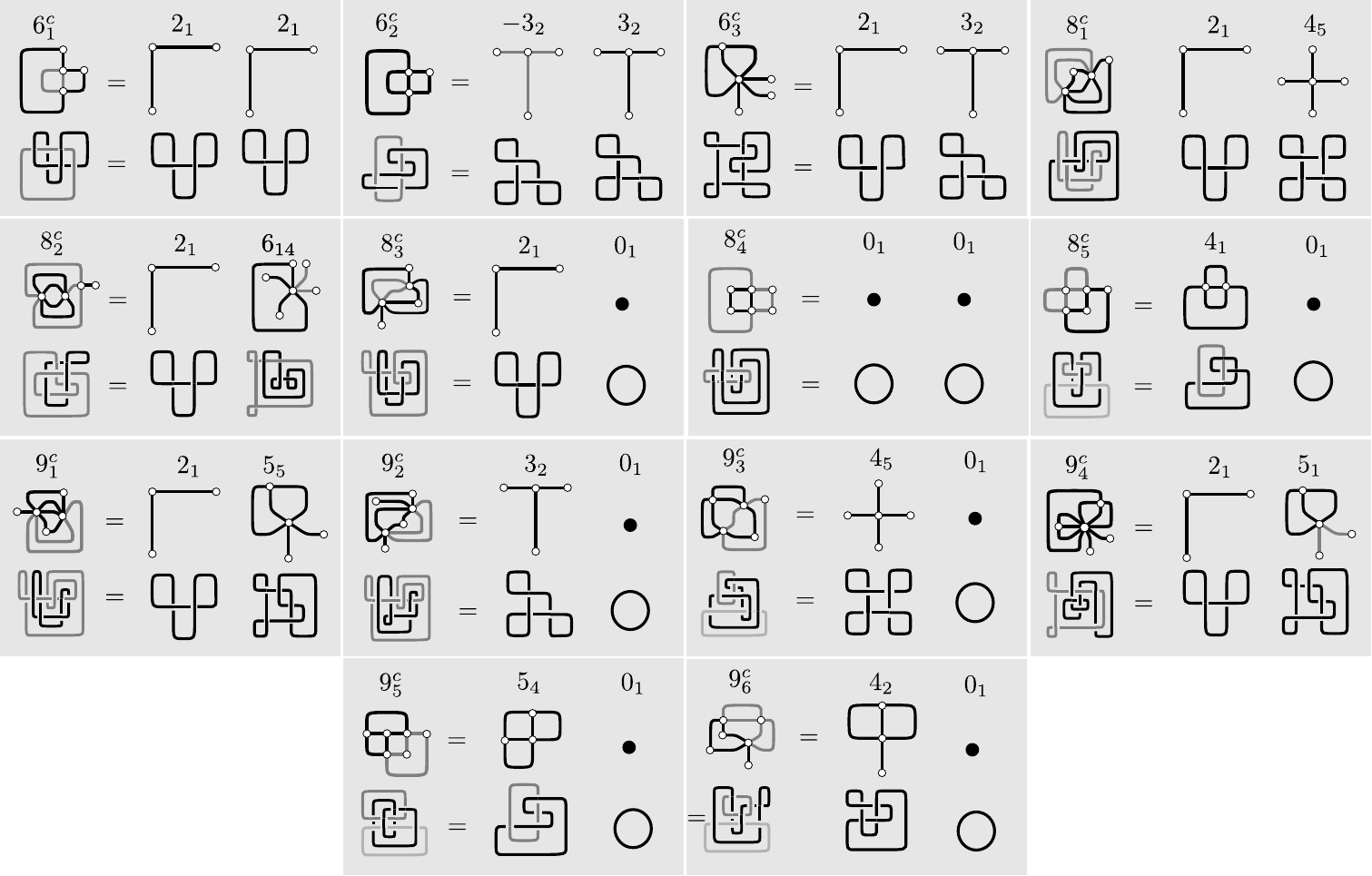}
\caption{\sf The 14 composite blink representatives up to 9 edges. 
Among the 3437 blinks in $U_9$ only the 14 $HG12QI$-classes (given above  by their 
representatives) induce
non-prime 3-manifolds. Detecion of non-primality is straighforward in  
gem (see Section \ref{sec:blinksgems}) language. 
It is shown by the existence of a {\em non-trivial combinatorial handle} in the gem.
This means a set of 4 distinclty colored edges so that any pair is in the 
same $ij$-gon, for the six possible choices of $\{i,j\}$. The combinatorial 
handle can be separating, in the case that deleting the 4 edges breaks the gem 
into two subgraphs having both more than 1 vertex (if one side has a single vertex, then the 
handle is trivial and is useless); or it can be non-separating,
in the case that it corresponds to connected sum with $\mathbb{S}^2 \times \mathbb{S}^ 1$.
There exists a linear algorithm to transform a blink into its canonical gem, 
see bottom part of Fig. \ref{fig:12to8move}. 
}
\label{fig:composite}
\end{center}
\end{figure}

Between the 1920's and beginning of the 1990's there was a feeling that good invariants to
3-manifold did not exist. Psicologically we were unprepared for the breakthrough that 
Witten's work in the connections between Physics and the Jones polynomial
(\cite{witten1989quantum})  implied to 3-manifold invariants. One of us, S. Lins,
in a sabatical leave to Chicado, lived intensely the excitement of these times and produced, 
in a join work with L. H. Kauffman \cite{kauffman1994tlr}, 
an effective completely combinatorial way for obtaining
the Witten-Reshetikhin-Turaev invariants, \cite{reshetikhin1991invariants}.  There are
at least four independent implementations, all agreeing, of the algorithms to get
WRT's. These invariants are the next ingredient after having a set of filtered blinks,
which misses no 3-manifold in certain classes.
This paper concludes the proof of the following theorem, 
\begin{theorem}
Let $\mathbb{B}$ be a connected blink with at most 9 edges and 
$\mathbb{M}^3$ be the closed, oriented, connected 3-manifold 
induced by $\mathbb{B}$.
There is a polynomial efficient algorithm which shows that
$\mathbb{M}^3$ is homeomorphic to exactly one of the 
3-manifolds induced by the 489  blinks below, or to 
one of the fourteen 3-manifolds induced by the non-prime 3-manifolds
of Fig. \ref{fig:composite}. 
Morever all of these are pairwise non-homeomorphic.
However being redundant, 
we also present the corresponding BFL's:
\newpage
\includegraphics[width=17.5cm]{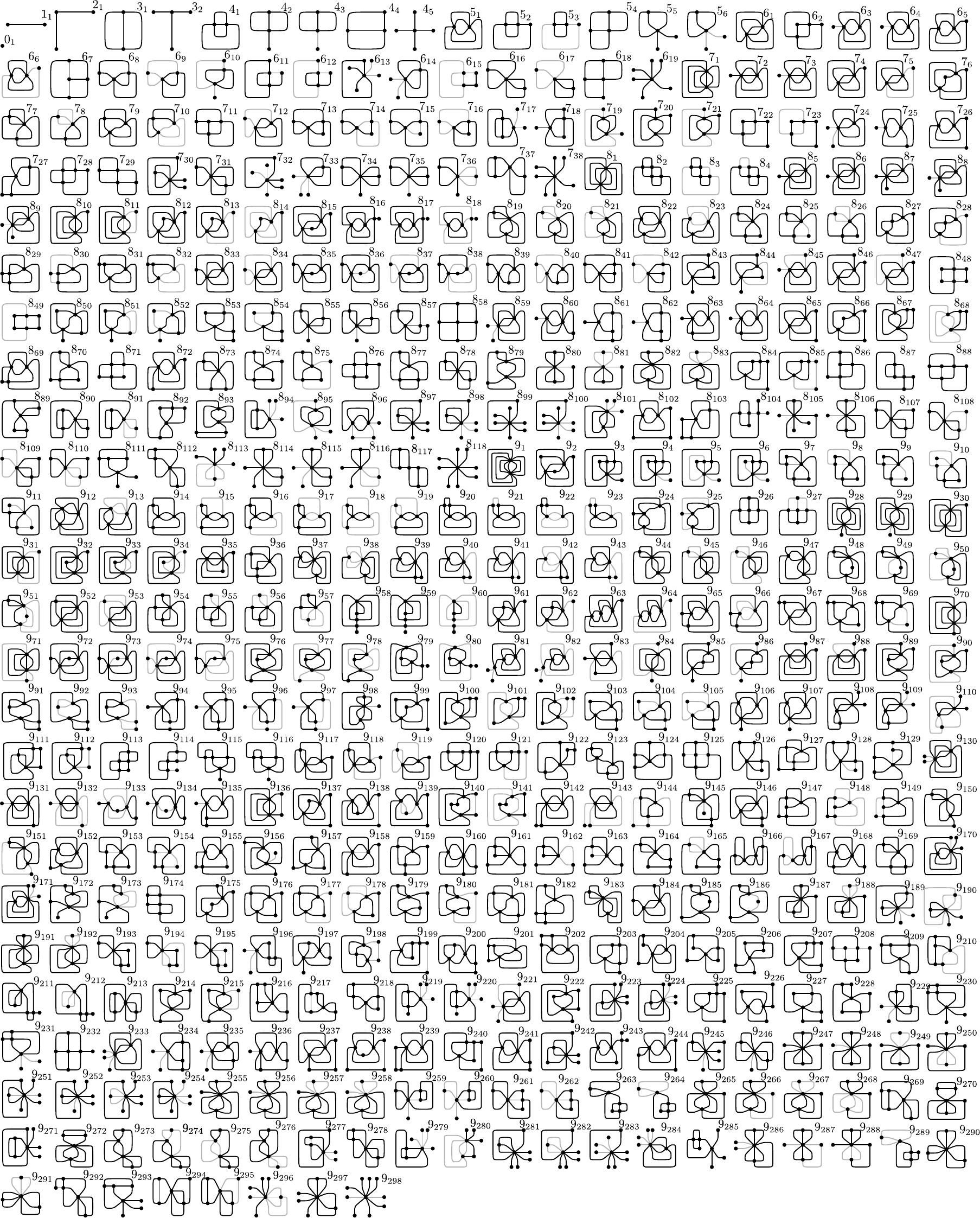}\\
\includegraphics[width=17.5cm]{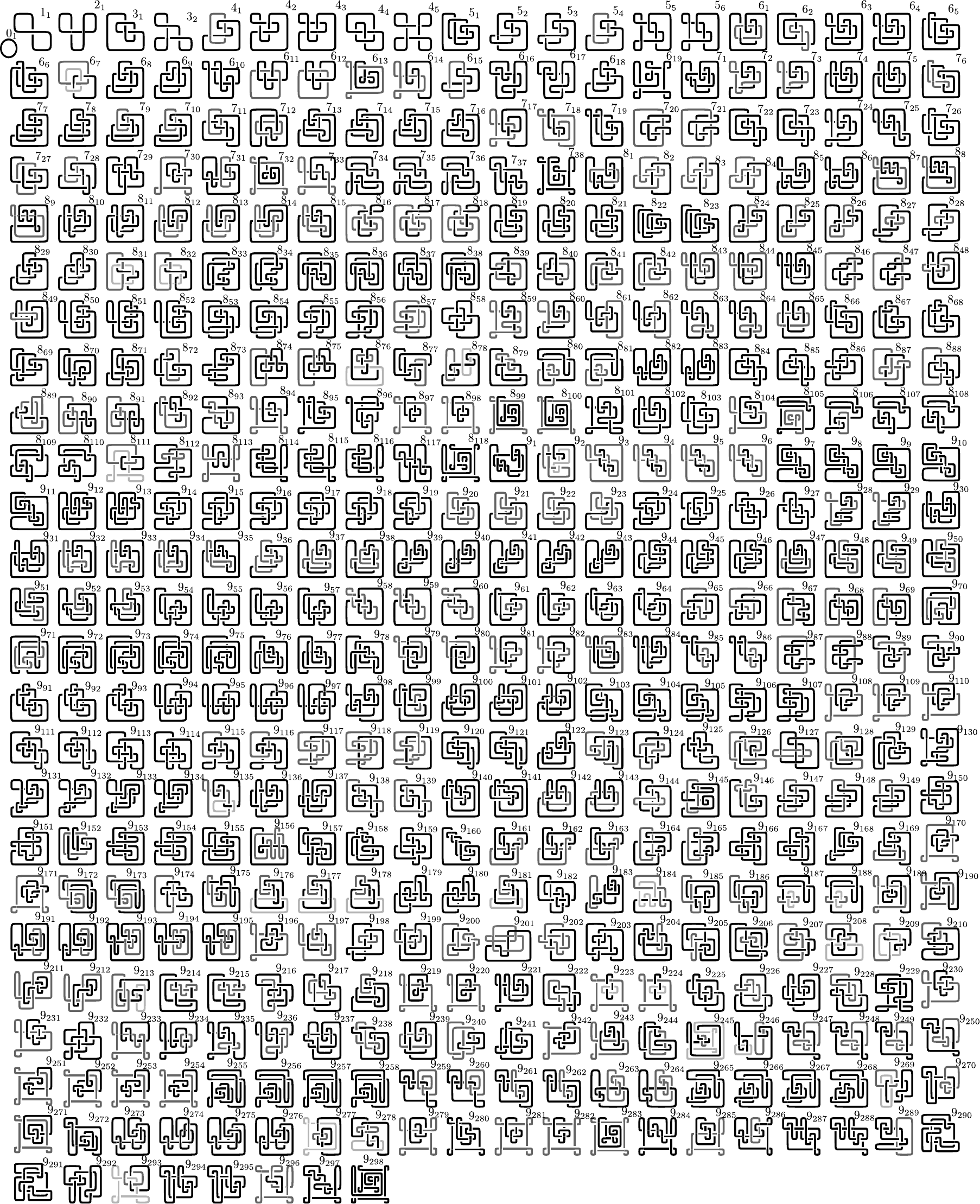} 
\label{theorem:census9}
\end{theorem}

In references \cite{kauffman1994tlr}, \cite{lins2007blink} and \cite{lins1995gca}
we have defined and show how a {\em blink}, that is, a plane graph with
an arbitrary bipartition of its edges (here presented as 
colors black and gray) induces a well
defined closed oriented 3-manifold. Moreover each such a manifold 
is induced by a blink (in fact, by infinite blinks).

\begin{figure}[H]
\begin{center}
\includegraphics[width=8.5cm]{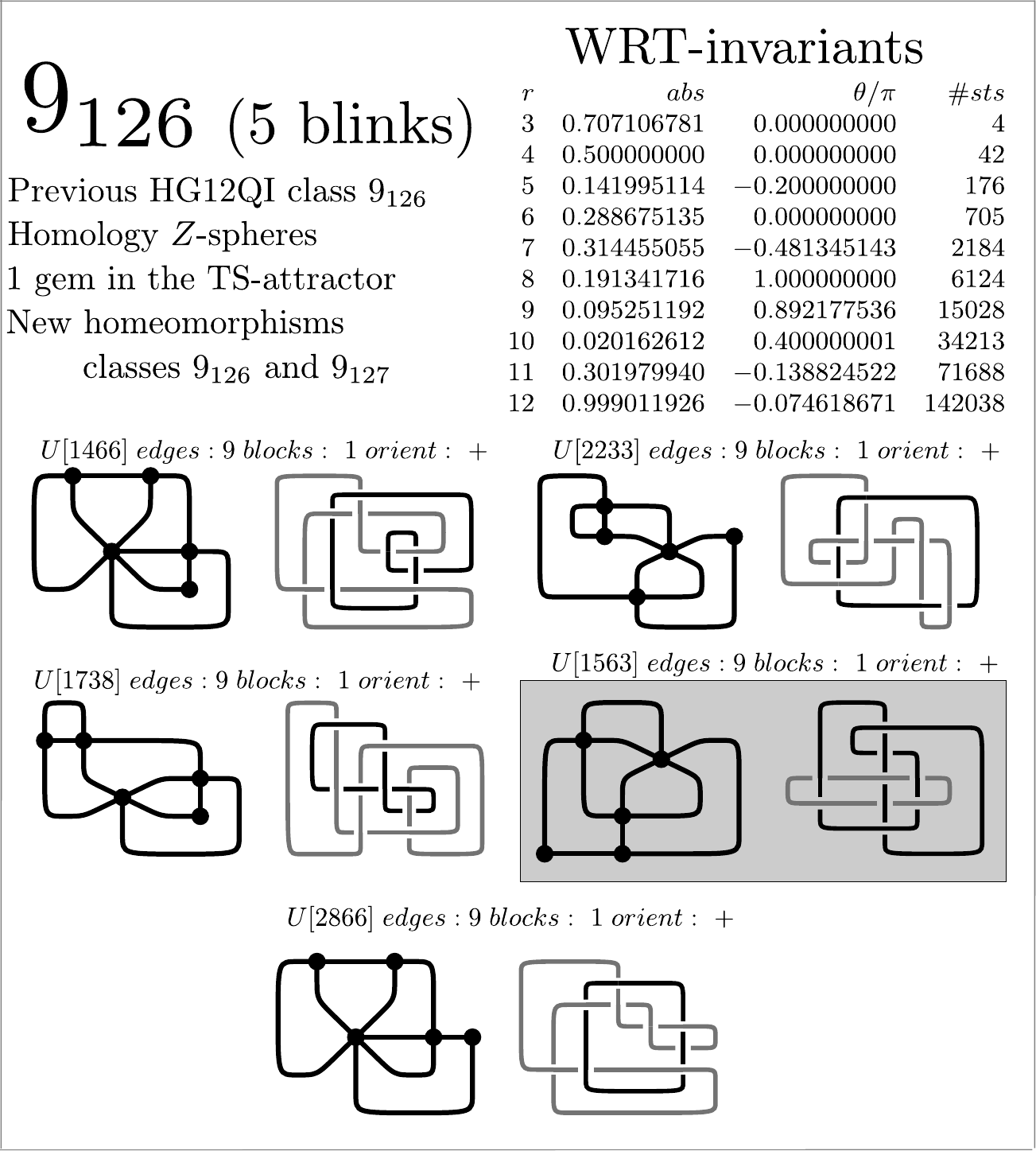}
\includegraphics[width=8.5cm]{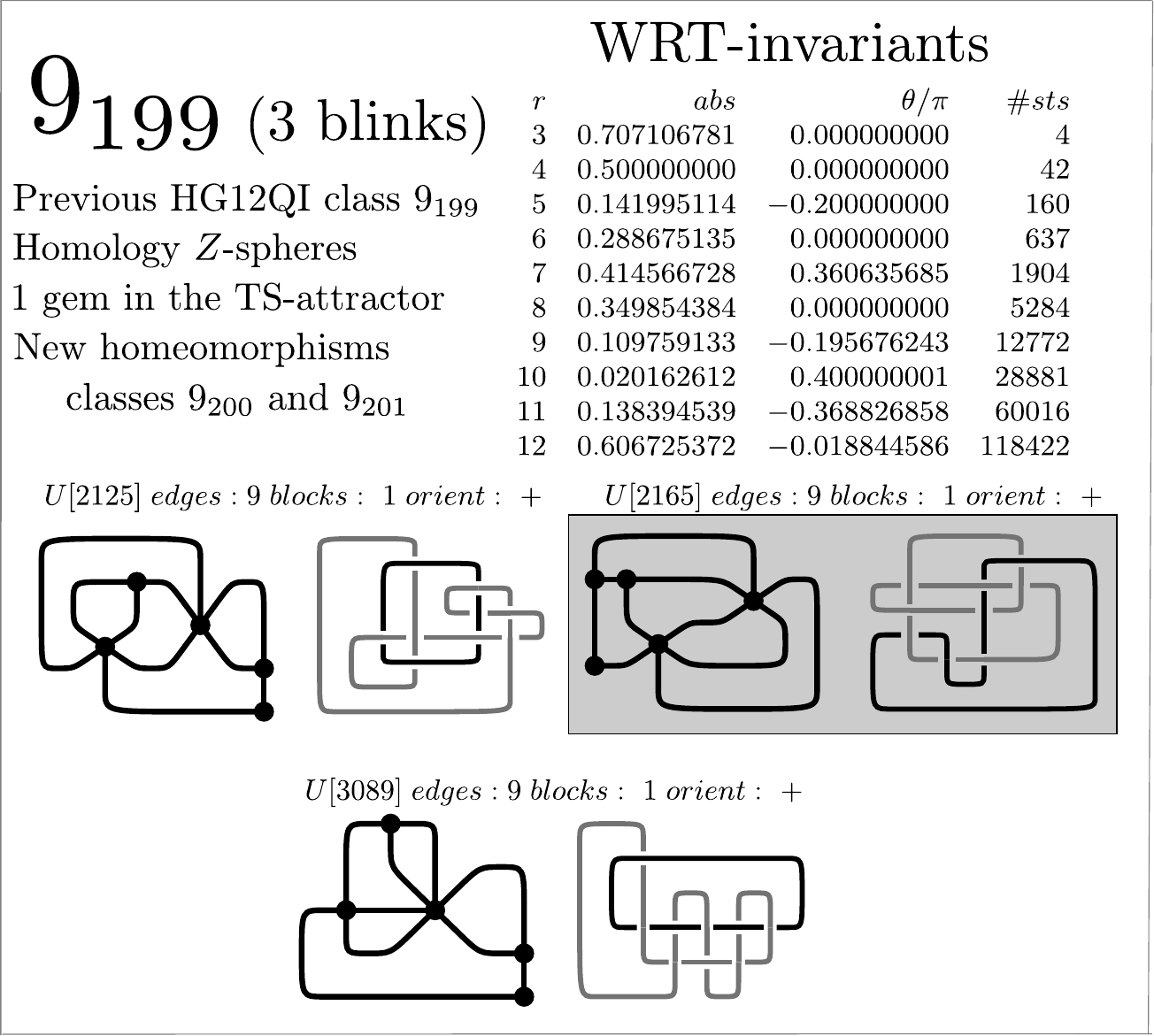}
\caption{\sf The two doubts left in L. Lins's thesis are solved:
in both cases the manifolds induced
by the shaded blink-link pairs are proved to be non-homeomorphic 
to the other in the same class, by geometric means. All the other pairs are proved 
to be homeomorphic by BLINK which keeps each such homeomorphism as a coded sequence which
is stored in a data basis, and, in principle reproducible at will. 
This shows that BLINK finds all available homeomorphic pairs and, in conjunction with
the length of the smallest geodesic (see the end of this section), provides
the topological classification of the 9-small 3-manifolds. Now we know that 
$U[1466] \equiv U[1738] \equiv U[2233] \equiv[2866]$ as they induce homeomorphic 3-manifolds
as proved by BLINK and by SnapPy. As for the shaded $U[1563]$ SnapPy/Gap/Sage proved that
its induced 3-manifold is non-homeomorphic. Similarly, $U[2125] \equiv U[3089]$ as proved
by Blink and SnapPy. As for the shaded $U[2165]$  SnapPy/Gap/Sage proved that the induced
3-manifold is non-homeomorphic. The $HG12QI$-class $9_{126}$ breaks into the topological classes 
$9_{126}$ and $9_{127}$. The $HG12QI$-class $9_{199}$ breaks into the topological classes 
$9_{200}$ and $9_{201}$.
}
\label{fig:nine126nine199}
\end{center}
\end{figure}

\begin{proof} {\bf (Of the Theorem above)}
 The bulk of the proof follows from L. Lins' thesis under the supervision of S. Lins, 
 \cite{lins2007blink}.
 In this work a theory for blink generation (missing no closed orientable 3-manifolds)
 is provided 
 by pure combinatorics, lexicography and topological filtering of duplicates.
 This resulted in a set $U_9$ with 3437 blinks. This universal set of 9-small 3-manifold
 was partitioned by homology and WRT$_{12}$ (the Witten-Reshetikhin-Turaev invariant 
 computed, as in \cite{kauffman1994tlr} for $r=3,4,\ldots,12)$ into 487 classes, named 
 $HG12QI$-classes. This is achieved by
 explicitly obtaining homeomorphisms between any two blinks in the same $HG12QI$-class.
 Each such homeomorphism is coded into a sequence of gems (defined in the next section)
 that is frozen in the data basis (Mysql) of BLINK forever, and, 
 in principle can be reproduced at will. Exactly 75633 gems were used in the 
 classifing sequences encoding homeomorphisms between the 3-manifolds of $U_9$.
 What remained to be done was to decide the status of the two $HG12QI$-classes
 $9_{126}$ and $9_{199}$ depicted in Fig. \ref{fig:nine126nine199}. After 6 years
 we posted these doubts as a Challenge in the arXiv, \cite{linslins2013A}. 
 In a quick feedback, many Mathematicians got interested in the Challenge (in the order of 
 their comments in the blog)
 \begin{center}
 \begin{verbatim}
 http://ldtopology.wordpress.com/2013/04/23/when-are-two-hyperbolic-
 3-manifolds-homeomorphic/#postcomment.
 \end{verbatim}
 \end{center}
 Comments were by H. Wilton, N. Dunfield, S. Friedl, M. Culler, D. Huberman, J. Berge, 
 C. Doria. Also, N. Dunfield M. Culler and C. Hodgson contacted one of us (S.Lins) via e-mail.
 And we got solutions for the two 6 year old doubts: 
 to the benefit of BLINK both pairs of 3-manifolds are non-homeomorphic.
 M. Culler, N. Dunfield and C. Hodgson sent distinct proofs of this
 fact. Thus there are 489 classes of non-homeomorphic closed 9-small 3-manifolds. More details
 of the distinction are given in Section \ref{sec:resol}.
 \end{proof}

Relative to page 109 of \cite{lins2007blink} the blinks of 
Theorem \ref{theorem:census9}
have receive two additions, the representative blinks
$U[1563]$ and $U[2165]$. Also the previous HG12QI-class $6_{5}$ became 
the homemorphism class $0_1$
corresponding to $\mathbb{S}^2 \times \mathbb{S}^1$. We have decreased 
by 1 the numbering of the 
$HG12QI$-classes $6_{6}, 6_{7}, \ldots, 6_{20}$ which become the 
homeomorphisms classes  $6_{5}, 6_{7}, \ldots, 6_{19}$:
the HG12QI-classes $9_{126}$ and $9_{199}$ of \cite{lins2007blink} split
into two homeomorphism classes. 
We observe that the blinks are enlarged in the appendix, showing them together with the
corresponding blackboard framed links. The notation $n_i$ attached to each blink below,
is the name of its homeomorphism class, not merely its HG12QI-class, as in \cite{lins2007blink}. 

\begin{figure}[H]
\begin{center}
\includegraphics[width=11.3cm] {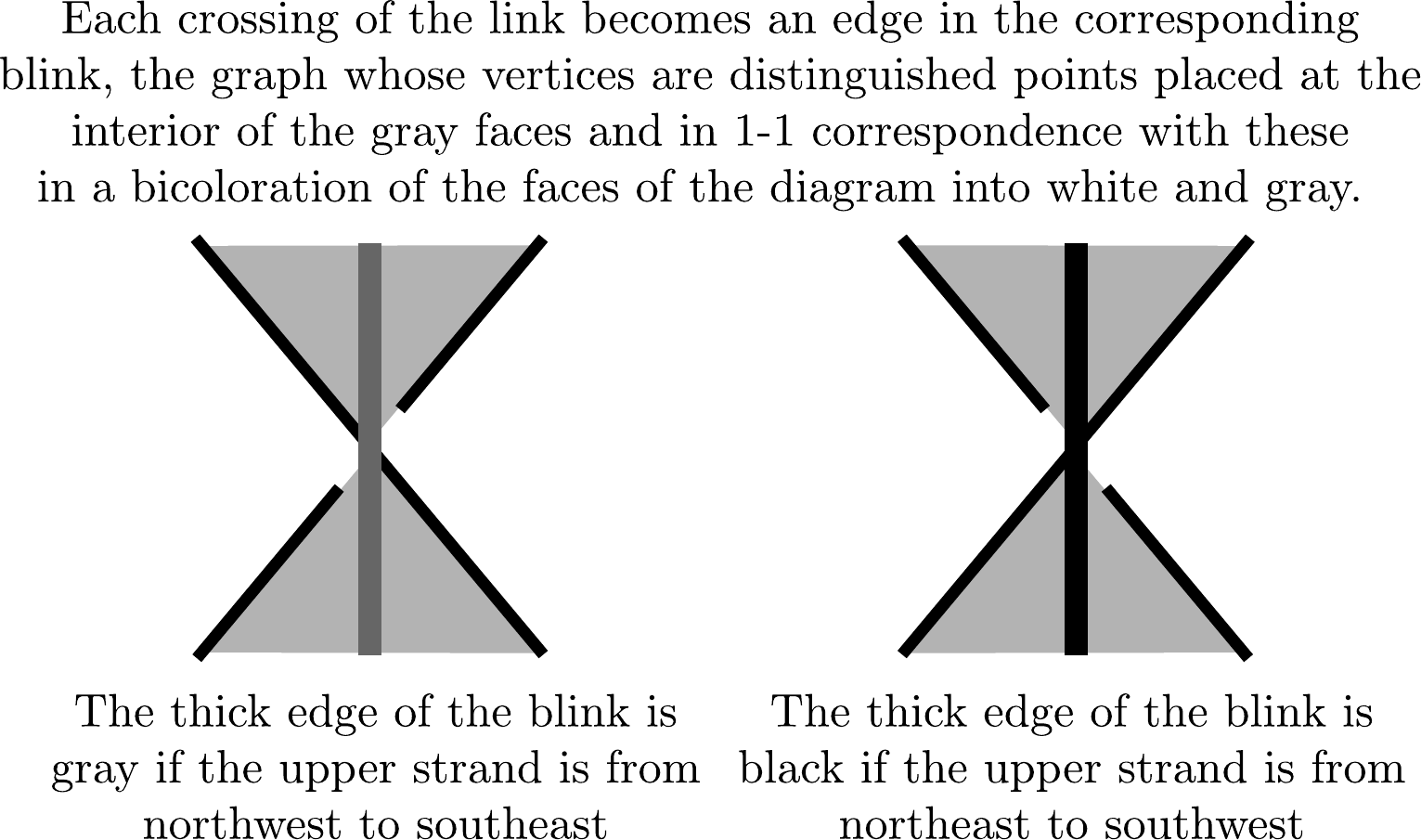}
\caption{\sf The linear algorithm to get a blink 
from a blackboard framed link. The inverse algorithm is also linear. 
}
\label{fig:linkblink}
\end{center}
\end{figure}

In Fig. \ref{fig:linkblink} we depict the function which associates a blink to a given
blackboard framed link via a checkerboard bicoloration of the faces of the link diagram.
The function has clearly an inverse and associated linear algorithms 
to go from BFL to link and from link to BFL are defined. 

\subsection{A word about the drawings in this work}
As important as the Mathematical content of this paper is the possibility of
easy visualization of the objects.
The elegant drawings of blinks and blackboard framed links produced by the software BLINK are 
possible due the groundbreaking algorithm of R. Tamasia \cite{tamassia1987egg}.
We could implement the drawings very fast because we had at hand the implementation of
network flow algorithms we had done for a project to solve {\em practical timetable (!) problems}.
This is an example of the unicity in Mathematics, advocated by L. Lovasz in his famous essay \cite{lovasz1998om}.
To get the drawings one has to apply three times the full strength of network flow theory,
\cite{ahuja1993network}.
The drawings BLINK present are in an integer grid and 
deterministically minimize the number of $\pi/2$-bents in the blackboarded framed links.
In particular, it permit us to deal with the unavoidable curls in the BFL's and the corresponding
loops and pendant edges which adjust the integer framings to self-writhes in
the best possible way: we do not care about them. 
The drawings for the companion blinks require a slight modification: it replaces each $p$-valent vertex $p>4$,
by a $p$-polygon inducing 3-valent ones. The final result is massaged a bit to
produce aesthetically pleasing and unambiguous figures in a homogeneous grid.

\subsection{Brief historical overview}

 It is amazing how much the picture has changed in the 80 years since the book by 
 Alexandroff and Hilbert.
 The progress initiated with the deep advances in the 1950's and 1960's, starting with
 the proof that 3-manifolds are triangulable by Moise (1952), \cite{moise1952affine}. Next 
 the presentation
 of them by framed links by Lickorish (1962) \cite{lickorish1962representation}.
 Following that Kirby presented its calculus for framed links (1978)\cite{kirby1978calculus}. 
 Starting in the early 1980's W. Thurston's provided gret 
 breakthroughs, developing his conceptual theory 
 on hyperbolic manifolds and of the geometrization conjecture. 
 In the final 1980's early 1990's Witten \cite{witten1989quantum} broke the psicological 
 barrier that there were no good invariants for 3-manifolds. Following that a number of 
 eastern European mathematicians like N. Reshetikhen, V. Turaev and O. Viro, 
 \cite{turaev1992state,reshetikhin1991invariants}
 using quantum groups were able to put in mathematical solid ground Witten's findings.
 One of us, S. Lins, was a witness of the excitement these developments caused. L. H. Kauffman
 and W. B. R. Lickorish discovered the relationship of the Temperley-Lieb algebra with 
 the new invariants,
 \cite{lickorish1991three}. Starting with a sabatical leave in to Chicago in 1990, S. Lins produced
 the joint monography with Kauffman \cite{kauffman1994tlr}, where blinks are first defined 
 and extensive WRT-invariant computations were obtained from the theory developed from scratch,
 independently and simpler than that of quantum groups.  
 The algorithms to obtain these invariants have been implemented by four independent reseachers,
 all agreeing.  In the early 
 2000's, G. Perelman revolutionized the field proving Poincaré's Conjecture and
 Thurston's Geometrization Conjecture. More recently in the 2010's, I. Agol is leading the field
 in this era post-Perelman. Of course, this is only a diagonal list of researchers. Many more have contributed and
 some are  extremely active in this era post-Perelman, \cite{friedl3}. 
 Currently there is a great amount of important reserch issues
 going on and these are exciting times for 3-manifold theory. See the recent essay of E. Klarreich
 in the Simons Foundation, \cite{Klarreich1012}.
         
\section{Blinks and gems}
\label{sec:blinksgems}
Unexplored simplicity. This was the reason for birth of this work many years ago.  
Repeating, a {\em blink} is a finite plane graph 
(that is, given embedded in the plane) together with an arbitrary bipartition
of its edges into black and gray. For completeness
we briefly recall the basic definitions of gem theory, leading to its 
definition, \cite{lins1995gca}
and to its calculus based in dipole moves, \cite{ferri1982crystallisation,lins2006blobs}.
A {\em 4-graph} $G$ is a finite bipartite 4-regular graph whose edges are partitioned into 4 colors,
0,1,2, and 3, 
so that at each vertex there is an edge of each color, 
a proper edge-coloration, \cite{bondy1976gta}.
For each $i \in \{0,1,2,3\}$, let $E_i$ denote the set of $i$-colored edges of $G$.
A $\{j,k\}$-residue in a $4$-graph $G$ is a connected component of the subgraph induced by $E_j \cup E_k$.
A 2-residue is a $\{j,k\}$-residue, for some distinct colors $j$ and $k$.
A {\em gem} is a 4-graph $G$ such that for each color $i$, $G\backslash E_i$ can be embedded in the plane 
such that the boundary of each face is a 2-residue. From a gem there exists a straightforward
algorithm to obtain a closed orientable 3-manifold, in two different, dual ways. 
Every such a manifold is obtainable in this way.
An unecessary big gem is obtained from a triangulation $T$ 
for a manifold by taking the dual of the 
barycentric subdivision of $T$. Here the colors correspond 
to the dimensions. 
Given a pair of vertices $\{u,v\}$ of a gem $G$ with $k$ linking edges $k \in \{0,1,2,3\}$, 
in $k$ distinct colors $K \subset \{0,1,2,3\}$ is called a {\em $K$-dipole} if $u$ and $v$ 
are in distinct $(\{0,1,2,3\}\backslash K)$-residues of $G$. A {\em dipole cancellation} is the 
operation that remove $u$, $v$ and their $k$ linking edges leaving $4-k$ pairs of pendant
edges which are reunited along the same colors. The {\em dipole creation} is the inverse move.
The dipole creations and cancellations do not change the induced 3-manifold.
It is possible to simplify substantially the gem corresponding to the dual 
of the barycentric subdivison of a triangulation by cancelling dipoles. 
These simplifications in that gem completely destroys the correspondence of 
colors and dimensions. Two gems induce the same 3-manifold if and only if they are linked 
by a finite sequence of moves where each term is either a dipole cancellation
or a diplole creation, 
\cite{ferri1982crystallisation,lins2006blobs}.

\begin{figure}[!h]
\begin{center}
\includegraphics[width=15cm] {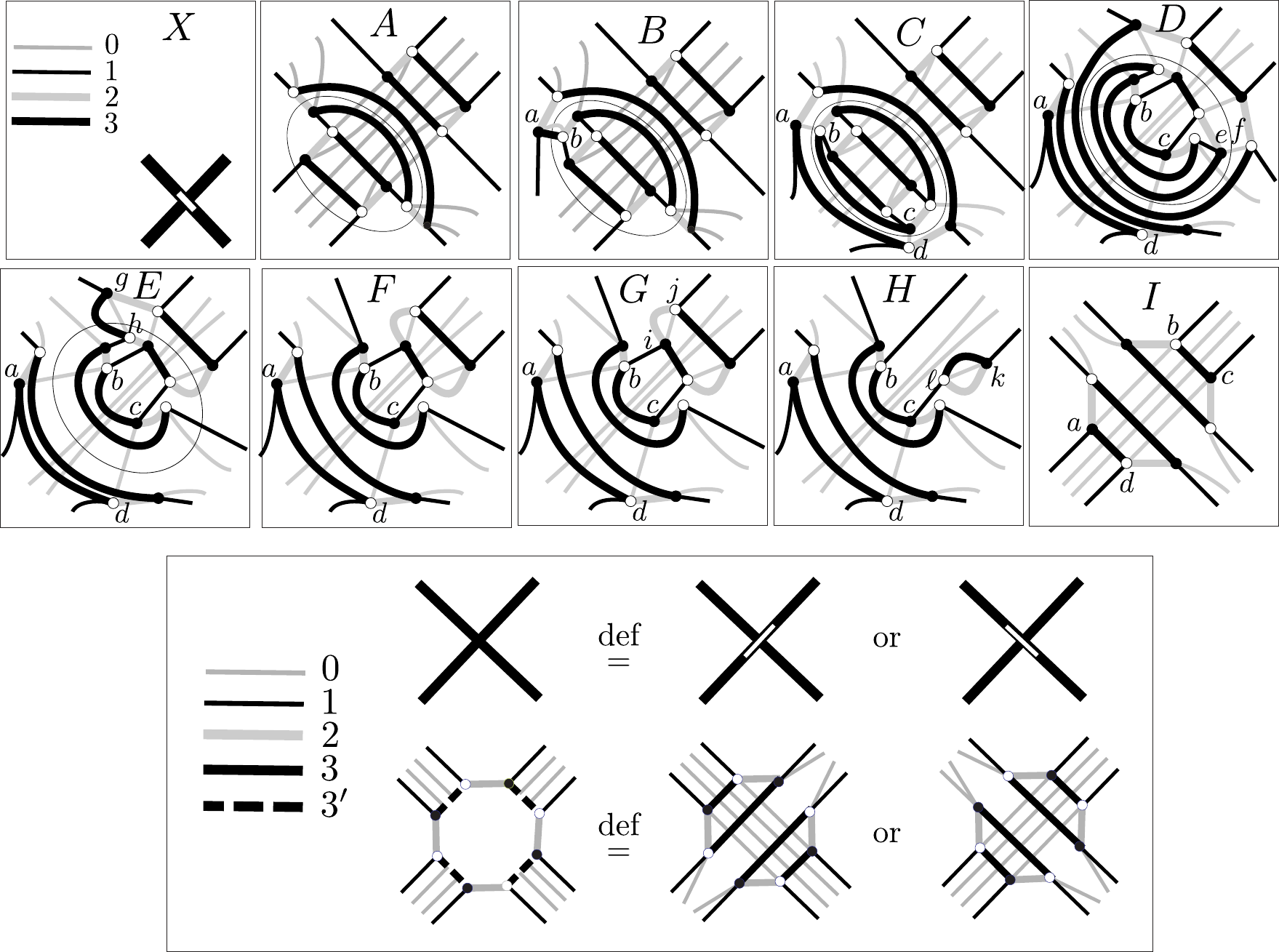}
\caption{\sf The passage $X \rightarrow A$ corresponds
to the $12X$-algorithm to go from a BFL to a gem, is defined and proved to be
correct in \cite{kauffman1994tlr}.
The $(12\to 8)$-simplification by dipole moves (passages $A \rightarrow B\rightarrow
\ldots \rightarrow I$) 
first appears in \cite{lins2007blink}. 
The final $8X$-algorithm to go from a BFL to a gem 
is depicted in the lower part of the figure. 
}
\label{fig:12to8move}
\end{center}
\end{figure} 

In our context, a very important property of gems is that, given any
blackboard framed link with $n$ crossings there is a linear algorithhm to 
present a gem having $8n$ vertices which induce the same 3-manifold.
This algorithm is displayed at the bottom part of Fig. \ref{fig:12to8move}.
This was proved originally in Kauffman-Lins monography, page 175 of 
\cite{kauffman1994tlr},
which provide a gem with $12n$ vertices. By using an specific $TS$-move of 
\cite{lins1995gca} a gem with $8n$ vertices is obtained, page 81 of \cite{lins2007blink}.
The $(12\to 8)$-simplification by means of dipoles is depicted 
the upper part of Fig. \ref{fig:12to8move}.

\subsection{A new $O(n^2)$-algorithm: from a resoluble gem to a blink}
\label{subsec:on2algo}

The present work is independent of Fig. \ref{fig:reducedblinkcalculus}, 
on Theorem \ref{theo:incredible} of Subsection{\label{subsec:surprise}}, on
Theorems \ref{theo:O2algo}, \ref{theo:framedlink} and on the Conjecture \ref{conj:remain}, 
of this Subsection, which can be skipped at no logical cost. 
But we recommend the reader not to do so,
because the Theorems and the Conjecture are stated to give a glympse and to clarify
our research in a broader appropriate context and timing.
Consider the following 9-element set of universal presentations of 3-manifolds:
\begin{center}
{\em $UP3M=\{$ 
(1) triangulation (E. Moise), 
(2) Hegaard diagram, 
(3) gem (Italian School of crystallizations \&S. Lins), 
(4) special spine (Matveev), 
(5) integer framed link  (Lickorish \& Kirby),
(6) blackboard framed link (L. Kauffman), 
(7) blink (S. Lins), 
(8) $\pm \infty$-fractional framed link (D. Rolfsen), 
(9) homology based framed link (SnapPy)
$\}.$}
\end{center}

The first 4 are {\em triangulation based types of presentations}. The last 5 are the
{\em framed link type based presentations}. Two triangulation based presentations 
are obtainable, one from the other, by inverse linear algorithms. If the integers are bounded 
in a presentation (5), then (5,6,7) are also linked by such algorithms. Presentations
(8,9) are more general than those of (5,6,7): the denominator of the fraction are
1 and the second number of the canonical homology basis are 1. Thus going from (5,6,7)
to (8,9), nothing needs to be done. The inverse algorithm, is easily to be shown to be finite,
although the complexity analysis is more difficult, \cite{rolfsen2003knots}. Fortunately,
these algorithms are not needed because we start with the blink, and so (5,6,7) are
equivalent up to linear algorithms.

It is straighforward to go from a blink to a gem by a linear algorithm: we take the
associated blackboard framed link and obtain the canonical gem by the linear algorithm of
Fig. \ref{fig:12to8move}. However, to go from a triangulated type presentation to a
framed link type presentation by a polynomial algorithm is, at present an open problem.

In a sequence of papers which S. Lins and R. Machado posed in the arXiv,
\cite{linsmachadoA2012,linsmachadoB2012,linsmachadoC2012} an $O(n^ 2)$-algorithm is
developed to produce a blink from a special kind of gems, named {\em resoluble}.
The importance of resoluble gems stems from the following theorem:
\begin{theorem}
\label{theo:O2algo}
There exists an $O(n^ 2)$-algorithm to go from a {\em resoluble} gem to 
a blink inducing the same 3-manifold.
\end{theorem}

The parts of the proof of this theorem are available at 
\cite{linsmachadoA2012,linsmachadoB2012,linsmachadoC2012}. Even though 
the work leading to the above result is still being polished, 
we feel that is important to mention it here because it implies that there is
an $O(n^2)$-algorithm to find the image of $\beta$ in Subsection \ref{subsec:surprise}. 
in the case that the 3-manifold is given by a
resoluble gem. The proper definition of resoluble gem is still unecessary complicated
in the posted papers. Here is the {\em currently adequate (we are still trying to improving it)} 
definition of resoluble gem. Let $G$ be a gem and
$u,v$ distinct vertices. For $j\in\{2,3\}$ we say that $\{u,v\}$ 
is a {\em $j$-double meeting, or a $dm_j$} if $u,v$ are in the same $\{0,j\}$-residue,
in the same $\{1,k\}$-residue and in distinct $\{2,3\}$-residues. 
Note that $\{0,1,j,k\}=\{0,1,2,3\}$.
Let $D_2$ be a set
of new edges (in color 4) linking $u$ to $v$, whenever $\{u,v\}$ is a $dm_2$ in $G$. 
Let $D_3$ be a set
of new edges (in color 5) linking $u$ to $v$, whenever $\{u,v\}$ is a $dm_3$ in $G$. 
Let $G_1=G$, $G_2$ be $G$ with its colors permuted under the 3-cycle $(1,2,3)$ and $G_3$ be
$G$ with its edges permuted under the 3-cycle $(1,3,2)$. We say that $G$ is a {\em resoluble gem} if
$G_i \backslash\{E_0 \cup E_1\} \cup (D_2 \cup D_3)$ is connected, for some $i \in \{1,2,3\}$.
Recently, in a still unpublished work, with the above adequate definition 
of resolubility the same authors have proved the following theorem:

\begin{theorem}
\label{theo:framedlink}
If a 3-manifold is given by a framed link then there exists a linear algorithm
to obtain a resoluble gem inducing the same manifold.
\end{theorem}

Two gems induce the same 3-manifold if and only if they are linked by a finite
sequence of {\em blob moves} and {\em valid flip moves}. This result is proved
for gems of arbitrary dimensions by S. Lins and M. Mulazzani in \cite{lins2006blobs}.
A {\em blob} $\{u,v\}$ is an $\{i,j,k\}$-dipole for 3 distinct colors $i, j, k$.
A {\em flip in a gem} is the recoupling of two equally colored edges so as to maintain 
the bipartiteness. Flips that transform the $I$-dipole $\{u,v\}$ into the 
the $I \cup \{j\}$-dipole ($j \notin I \subset \{0,1,2,3\}$) and their inverses
are named {\em valid flips}. What is needed to link by or $O(n)$ or $O(n^2)$  algorithms 
any two of the following set of universal presentations of the same 3-manifold $(UP3M)$,
is the following Conjecture:

\begin{conjecture}
\label{conj:remain}
 A resoluble gem remains so under a blob move, or a valid flip move and straightforward
 simplifications using $TS\rho$-moves, \cite{lins1995gca}.
\end{conjecture}

This Conjecture is one of those annoying results that we believe to be true
(there is too much freedom in the situation) but which stubbornly refuses to surrender.
We still did not put enough effort in trying to prove it, mainly because of lack of time.
For the present status of this whole theory what imports is that simplifying a gem by
TS-moves we invariably obtain a resoluble gem. For our computational purposes
there is an $O(n^ 2)$-algorithm to produce a blink inducing the same space as the one
given by an arbitrary member of $UP3M$. At the level of the catalogue, simplifying the 
gems under $TS_\rho$ produce resoluble gems. Note that, as a bonus
we gain the capability of computing the Witten-Reshetikhin-Turaev invariants
\cite{reshetikhin1991invariants, kauffman1994tlr}
from all the members of $UP3M$, what is currently impossible, because the 
WRT-invariants are only computable from framed link type presentions (6,7).

We stress our belief that many interesting and deep consequences to be investigated
of Theorem \ref{theo:incredible} 
to the topology of 3-manifolds and to the combinatorics 
of plane graphs ought to exist. Here is an example: what does it means,
in terms of blinks a closed orientable 3-manifold having 
a homogeneous Riemannian geometry?

 \subsection{Blackboard framed links}
 Some pictures and some definitions of this subsection are reproduced 
 from Section 12.1 of \cite{kauffman1994tlr}.
 They are included here for completeness since the concept is central in this paper.
 Instead of working with the link in $\mathbb{R}^3$ it is to our 
 advantage to work with a {\em general position decorated projection} of the link
 into $\mathbb{R}^2$. Such projections are named {\em link diagrams}. 
 General position means that the pre-image of any
 point of $\mathbb{R}^2$ in the link is at most two points. A point that has
 two pre-image is a transversal crossing of two strands in the projection. Decorated means
 that we keep the information of which strand is the lower and which is the upper,
 usually by removing a small segment of the lower strand centered at the crossing.
 
 If $L$ is a link in $\mathbb{R}^3$ a {\em framed link based at $L$} is an embedding 
 $f: \mathbb{S}^1 \times [0,1] \longrightarrow \mathbb{R}^3$
 so that $L=f(\mathbb{S}^1 \times \{0\}).$ {\em Regular isotopy} is the class of
 link diagrams up to Reidemeister moves 2 and 3. Framed links are represented 
 by regular isotopy classes of link diagrams, with the extra move, 
 the {\em ribbon equivalence}, which
 is generated by the relations 
 \raisebox{-4mm}{\includegraphics[width=6.7cm]{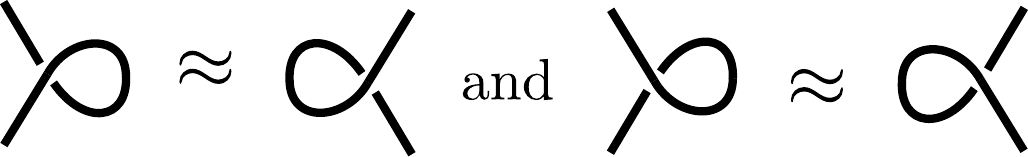}} .
 
\begin{figure}[H]
\begin{center}
\includegraphics[width=11cm]{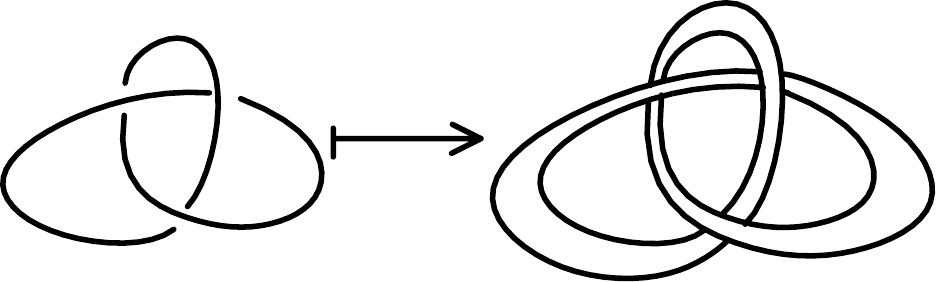} 
\caption{\sf How the tips of the normal vector field draws a second parallel component}
\label{fig:doubledtrefoil}
\end{center}
\end{figure}

 With these extra relations, the ribbon equivalence classes of link diagrams 
 represent the framed links via the {\em blackboard framing}. This framing is 
 obtained by placing a vector field normal to each link component 
 with all the normals in the plane indicated by the diagram. We can indicate
 this this framing by sketching the tips of the normals as a second component,
 so that each component is indicated by an immersed band,
 see Fig. \ref{fig:doubledtrefoil}.
   
The reason for the ribbon equivalence is manifest in that the framings corresponding to these
two curls (with distinct winding number in the plane) are ambient isotopic:
$$\includegraphics[width=11cm]{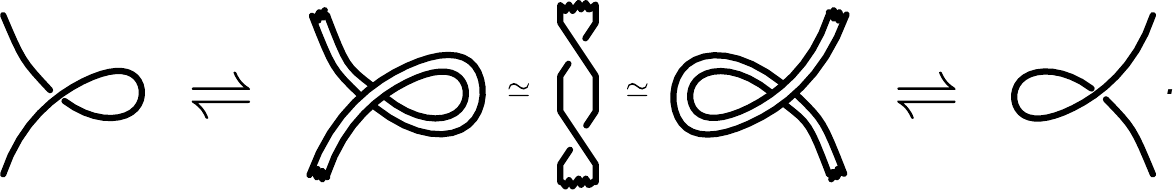}$$ 
We finish this subsection with a representative example of how to construct the closed 3-manifold
from a $(p,q)$-parametrized framed link (Dehn's $(p,q)$-filling). 
Suppose that the link diagram $K$ has just one component
and 3 curls corresponding to framing number 3 in $\mathbb{R}^3$.
Then $\mathbb{M}^3(K)$ is obtained by sewing 
$m \subset \mathbb{D}^2 \times \mathbb{S}^1$ to $\pm p.\theta_k+q.\lambda_K \subset
\mathbb{S}^3 \backslash (\mathbb{S}^1 \times \mathbb{D}^2)^{\tiny{o}},$ 
where $X^{\tiny{o}}$ means the 
interior of $X$. In the particular case of the blackboard framed link,
we must have $p=3$ and $q=1$. In this case the sewing, 
forms the lens space $L_{3,1}$. 

 \begin{figure}[H]
\begin{center}
\includegraphics[width=13cm]{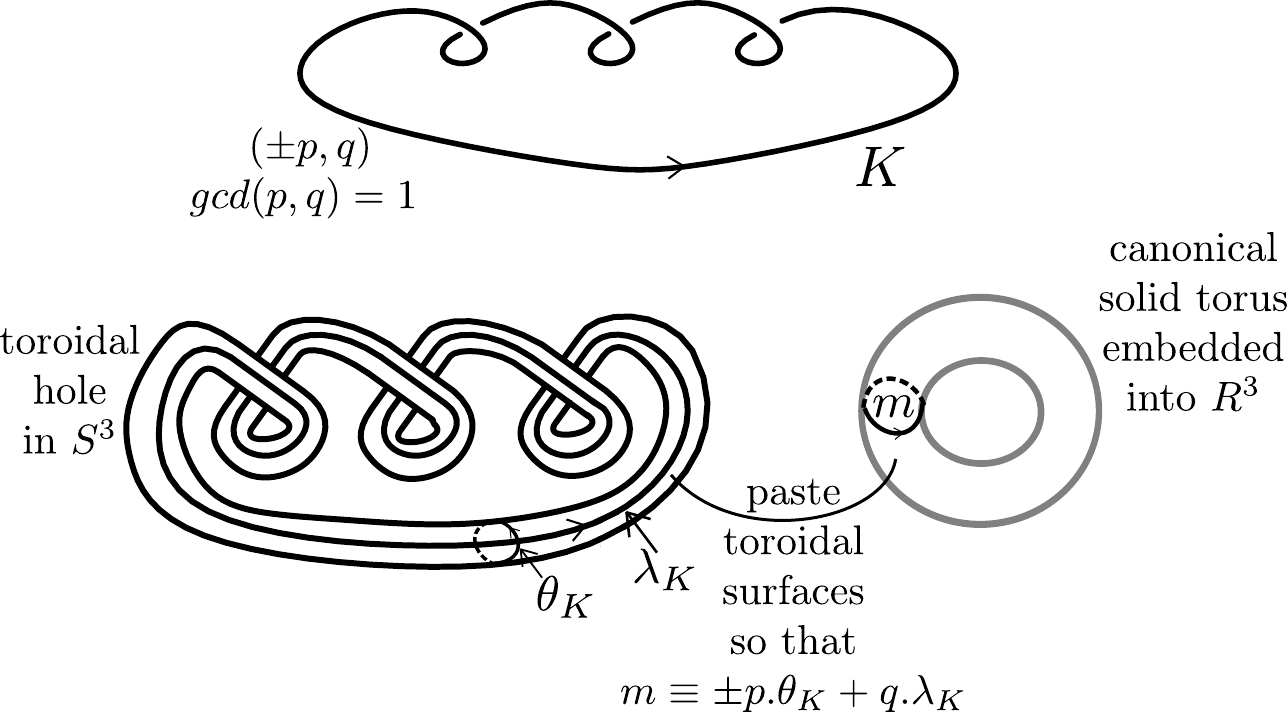} 
\caption{\sf Note that the longitudinal curve $\lambda_K$ in the boundary of the
toroidal hole never touches the boundary of the projection of 
the band formed by $K$ and its copy. Any $\epsilon$-neighborhood of the band removed form
a toroidal hole. 
The original knot component $K$ has linking number -1 with $\theta_K$.
Surgeries with $q=1$ are sufficient to generate
any closed oriented 3-manifold. Those are the surgeries involved in Lickorish's
original paper \cite{lickorish1962representation}. An algorithm to update the link 
so as to replace a non-negative 
general integer $q$ in each component by $q=1$ (while maintaining the induced 3-manifold) 
is given in Rolfsen \cite{rolfsen2003knots}. Note that in the case of blackboard framed links
the specification of $(\pm p,q)$ becomes redundant: the decorated projection of the 
link is all that matters.
}
\label{fig:General_pq_DehnFilling}
\end{center}
\end{figure}

\subsection{The complementary roles of blinks and gems}

\begin{figure}[H]
\begin{center}
\includegraphics[width=11cm]{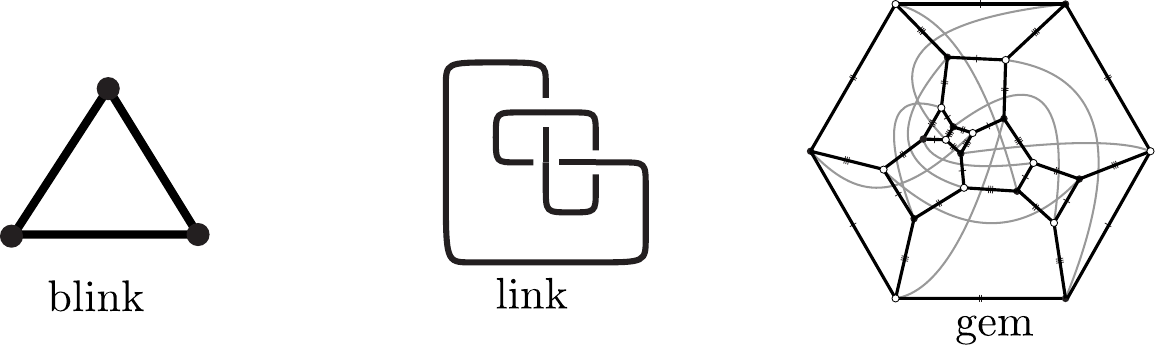} 
\caption{\sf The minimum blink, minimum link and minimum gem inducing the 
binary tetrahedral space}
\label{fig:blinklinkgem}
\end{center}
\end{figure}

 Blinks are in 1-1 correspondence with blackboard framed links which in turn encodes every 
 closed oriented 3-manifold. If we consider these three 
 encodings of the same 3-manifold given in Fig. \ref{fig:blinklinkgem}, 
 the blink is the one that has the smallest ``perceptual complexity''. 
 The common manifold of this example is the binary tetrahedral space defined
 by $\mathbb{S}^3$ (as a topological continuos group) 
 up to the action of the non-comutative binary 
 tetrahedral group $\langle 3,3,2\rangle$, 
 (see \cite{magnus1966combinatorial}) which has 24 elements. The 3-manifold
 $\mathbb{S}^3/\langle 3,3,2\rangle$ has an spherical geometry. Its attractor
 (definition given shortly) consists
 of the single gem depicted at the right of Fig. \ref{fig:blinklinkgem}.

 Blinks are very easy to generate recursively. They have a rich
 simpifying theory which permits the generation of 3-manifold catalogues. 
 Moreover, their isomorphism problem is computationally  simple. 
 Being blinks so good, why do we need gems? The 
 answer is that to prove that blinks induce the same 
 3-manifold (when they do) remains (probably even with the new coin calculus)
 a very difficult problem. It is 
 straighforward to obtain a canonical gem from a blink, \cite{kauffman1994tlr,lins2007blink}.
 Proving that two gems induce the same
 3-manifold (when they do) is much easier because of the rich simplification 
 theory of them (based
 on at least 4 intertwined planarities) leading to the {\em attractors 
 of 3-manifolds}, \cite{lins1995gca}. These are defined as the set of gems inducing 
 the manifold having the minimum set of vertices. The important computational 
 point is that the attractor
 usually have few members, particularly if the manifold has a uniform geometry
 euclidean, spherical, hyperbolic. 
 
  Blinks are very good at proving that two manifolds are distinct, because 
 from them we can extract
 the WRT-invariant, which is a very strong invariant, yet not a complete one. 
 These invariants currently 
 do not have a direct computation neither from triangulations, neither from gems,
 neither from Heegaard diagrams, neither from special spines of Matveev.
 In this respect see subsection \ref{subsec:surprise}, where we report progress
 on this issue.

 Gems and  blinks collaborate in a symbiotic dance to decide 
 (at a computational level) whether two 3-manifolds are or are not homeomorphic. 
 And many types of census  become available!  
 We present here our contribution to the topic.
 It is placed in the confluency of two deep passions of the authors:
 the study of closed orientable 3-manifolds and the study of plane 
 graphs. Here, we mean to provide a strategy for a 
 segmented answer of Hempel's questions posed at the introduction. 
 The closed oriented 3-manifolds are partitioned
 by the number of edges in a minimum encoding of them by a blink.

\begin{figure}[H]
\begin{center}
\includegraphics[width=6.5cm]{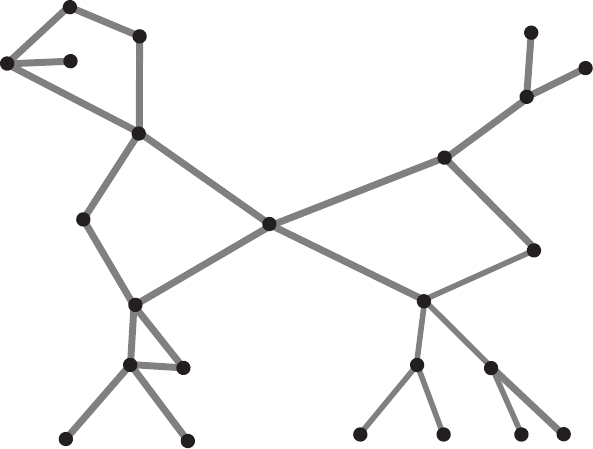} 
\raisebox{1.7cm}{
$\mapright{unique\  representative\ algorithm} {yielding\ the\ canonical \ form \ of \ the\ blink}$ }
\includegraphics[width=3cm]{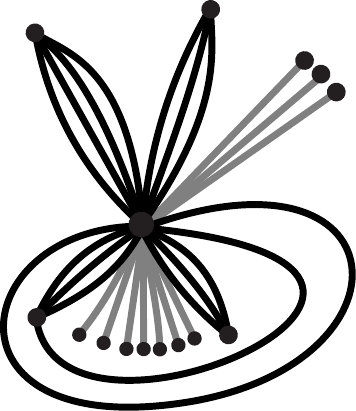} 
\caption{\sf Blink unique representative obtained by an $O(n\log n)$-algorithm 
implemented by BLINK: monochromatic 
doglike blink with 27 edges 
and its representative with 25 edges.}
\label{fig:doglikeblink}
\end{center}
\end{figure}

Any plane drawing whatsoever (see Fig. \ref{fig:doglikeblink})
of a graph with an arbitrary bipartition of its edge set, that is, a blink,
corresponds to a unique closed oriented 3-manifold via the associated 
blackboard framed links. An important aspect
about blinks is that each one possesses an easily obtainable {\em canonical form}
inducing the same 3-manifold: it is named the {\em representative of the blink} and is obtained
by lexicography from a small number of conventions, fixed in advance. 
This is explained, with a great amount of details, 
in L. Lins' thesis, \cite{lins2007blink}. 
 
 Lexicography is used to define a representative unique plane graph (a canonical form) 
 for each closed oriented 3-manifold. We explicitly solve the segmented problem 
 up to 9 edges, see Theorem \ref{theorem:census9}. This work provides 
 an efficient algorithm to make available the canonical form of any closed orientable 3-manifold
 induced by plane graphs at the current level of the catalogue (currently 9 edges) and, theoretically, this
 could be extended 10, 11,\ldots, $n$, for arbitrarily large $n$. Our theory gives a road to effectively
 name each 3-manifold classified by some set of invariants $INV$. We have a universal set of object, the blinks up
 to $n$ edges, which can be partitioned by these invariants. The $INV$-classes are then tried to be broken into
 homeomorphisms classes. New invariants are then discovered and added to $INV$ making them homeomorphisms classes 
 of $n$-small manifolds. The difficult cases are going to appear naturally and they lead to enhancement of
 the theory. It is not at all impossible that this process stops and we get $INV$ so 
 that the $INV$-classes be proved to  be homeomorphisms classes for all $n$. The point we want to make is that 
 good examples (hard to find, here exemplified by the $HG12QI$ classes 
 $9_{126}$ and $9_{199}$) are important in obtaining progress in a general theory.
 A classifying $INV$ for the 9-small 3-manifolds is
 $INV=\{{\ homology}, W\hspace{-0.5mm}RT_{12}, { \ length \ of\ smallest\  geodesic}\ \}.$

\section{Generating the $k$-universal set of blinks $U_k$}
\label{sec:generation}

A set of blinks is said to be $k$-universal if its members induce all $k$-small prime
3-manifolds. Note that the components of a disconnected blink correspond to the summands 
in the connected sum of the 3-manifolds induced by the components of the blink. 
Therefore, our minimal sets of $k$-universal blinks are
composed of connected blinks. It is not difficult to implement an algorithm 
that produces a specific
$k$-universal set $U_k$. We start by (1) constructing all the 2-connected graphs having up to
$k$ edges, forming a set $A$. This is a standard procedure in graph theory.
Then we take (2) all unions of blocks 
having a common vertex so that their edges 
set has at most $k$ edges. Regular isotopy shows that we can restrict to having all
blocks with a common vertex. And, by lexicography we do this in a unique way.
Denote by $B$ the set of these combined blocks.
Next (3) we consider all bicolorations of the edges so that
the number of black edges is not smaller than the number of gray edges. We put the blink 
in the third set only if it coincides with its representative. This forms the set $C$.
Finally we use some topological filtering to decrease the size of $C$. The surviving blinks 
forms the set $D$, which is our $U_k$. A final comment we make about the generation is
that the connected blinks are implemented as embedded in $\mathbb{S}^2$. The external face
of a blink is well defined: namely is the face indexed with number 1 in the 
numbering which produces its code. At any rate, changing the external face
is easily done in the coin calculus: it involve regular isotopy moves, 
one ribbon move and
one Whitney trick simplification.

\begin{figure}[H]
   \begin{center}
      \leavevmode
      \includegraphics[width=17cm]{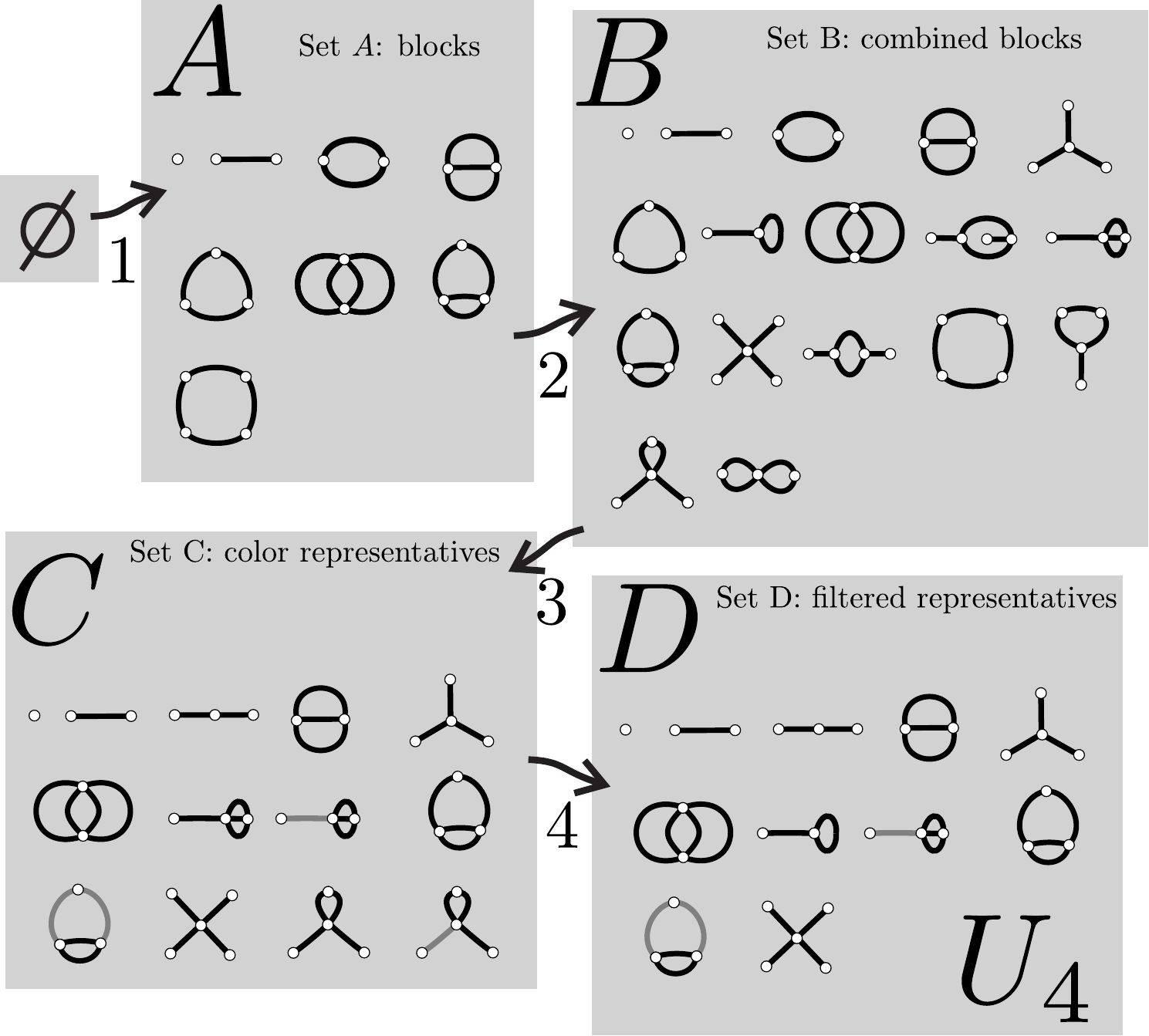}
   \end{center}
   \vspace{-0.7cm}
   \caption{Pipeline for obtaining, in 4 steps,
   the set $U_k$ exemplified at $k=4$. Note that, starting
   with the empty set, the input of a phase is the output of the previous phase, that it,
   we have a pipeline. We must
   include the single vertex, as it corresponds to the unknot, inducing 
   $\mathbb{S}^2 \times \mathbb{S}^1$. The cardinality of $U_4$ is only 11.
   BLINK produced, in its data base, the sets $U_9$ and $U_{10}$ with respectively
   3437 and 17948 blinks.
   }
   \label{fig:4pipeline}
\end{figure}

\begin{figure}[H]
\begin{center}
\includegraphics[width=14cm] {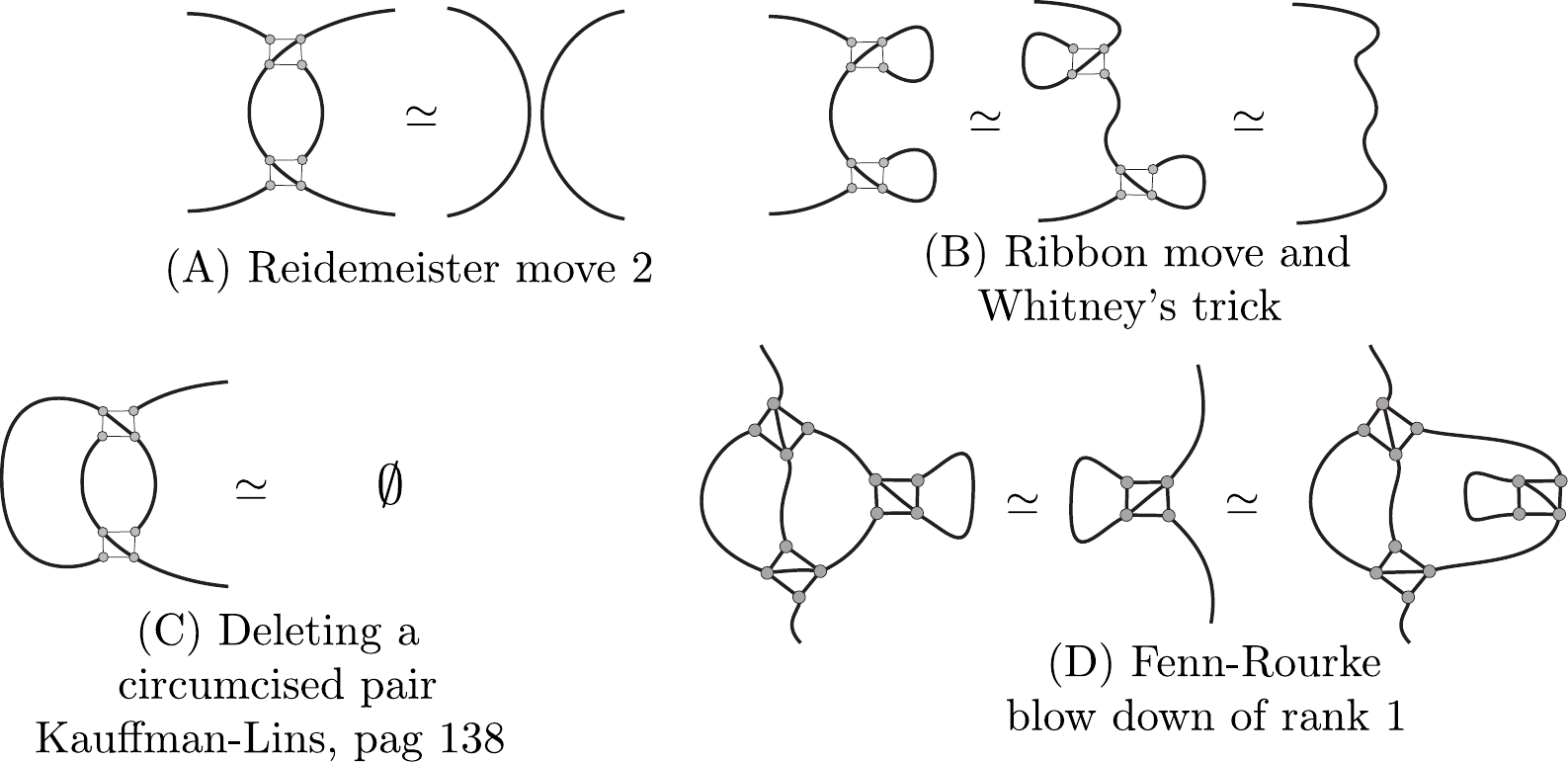}
\caption{\sf The topological filtering (at the level of blackborad framed links)
in the generation of $U_k$. A final ingredient that we use in the topological 
filtering is to drop blinks that produce non-prime 3-manifolds. They can be
detected easily by finding a separating or non-separating handle 
in the associated gems: just a 
set of 4 differently colored edges pairwise in the same $ij$-gon,
for the six types of such gons.
There are only 14/3437 composite 
$HG12QI$-classes of blinks in $U_9$.
}
\label{fig:SomeTopologicalFiltering}
\end{center}
\end{figure}

The topological filters that we use are depicted in Fig. 
\ref{fig:SomeTopologicalFiltering}.

\section{Unveiling the mystery of the two doubts in L. Lins thesis}
\label{sec:resol}

The topological classification of the 9-small spaces was nearly completed in 
\cite{lins2007blink}. This work develops a theory for generating a distinguished
set of blinks named $U_n$ and indexed lexicographically, $U_n[i]$ is the $i$-th such blink.
This has been revewed
The relevance of $U_n$ is that it misses no closed, orientable, prime and irreducible
3-manifold which is induced by a blink up to $n$ edges. 

The 3-manifolds of \cite{lins2007blink} are classified by homology and the quantum WRT$_r$-invariants
$r=3,\ldots,u$, with $10$ significant decimal digits forming $HGuQI$-classes of blinks. 
Our algorithm for computing the 
$WRT_r^u$-invariants are based on the theory developed in \cite{kauffman1994tlr}.
After 6 years we have put our doubts as a Challenge to topologists and 
group algebraists, \cite{linslins2013A}. They were quickly solved by some researchers among others
M. Culler, N. Dunfield and C. Hodgson, at least in two different ways. They proved that the two pairs of
manifolds which were left unresolved by BLINK are indeed non-homeomorphic. All the other pairs
in the $HG12QI$-classes $9_{126}$ and $9_{199}$ have been checked to be homeomorphic, as BLINK proved
6 years ago. The new solutions were obtained using the software \cite{snappy}, which uses 
the kernel of \cite{weeks2001snappea}. They also use GAP (\cite{gap2002gap}) and Sage (\cite{sage2012}). 
 
The first solution we obtained was by C. Hodgson using length spectra techniques, 
based in his
joint paper with J. Weeks entitled {\em Symmetries, 
isometries and length spectra of closed hyperbolic three-manifolds} 
(\cite{hodgson1994symmetries}). By using SnapPy Hodgson showed that even though the 
quantum WRT-invariants as well as the volumes of the 
hyperbolic $Z$-homology spheres induced by the blinks
$U[1466]$ and $U[1563]$ are the same, the {\em length of the 
smallest geodesics of them are distinct}. As for the other pair of blinks, $U[2125]$ and $U[2165]$, 
the same facts apply. Here is a summary of Hodgson's findings extracted from the SnapPy session
that he kindly sent us. As he explains: {\em ``The output of the length spectrum 
command shows the {\em complex lengths}
of closed geodesics --- the real part is the actual length and the imaginary 
part is the rotation angle
as you go once around the geodesic.''}\\
\begin{center}
\center{Class $9_{126}$:}
\begin{verbatim}
First geodesic of U[1466]: 1.0152103824828331+0.39992347315914334i.
First geodesic of U[1563]: 0.9359206605025168+2.333526236965665i.
Volume of both manifolds: 7.36429600733.
\end{verbatim}
\center{Class $9_{199}$:}
\begin{verbatim}
First geodesic of U[2125]:  0.8939075859248593+0.761197185679321i.
First geodesic of U[2165]:  0.7978548001747316+2.9487425029345973i.
Volume of both manifolds: 7.12868652133.
\end{verbatim}
\end{center}
Both N. Dunfield and M. Culler also classified topologically all the 3-manifolds
in the two $HG12QI$-classes of 3-manifolds,
by $k$-coverings techniques, an independent check of
Hodgson's proof. We thank to the these 3 researchers: M. Culler, N. Dunfield and C. Hodgson
for unveiling our 6-year old mystery: as we believed, both  $HG12QI$-classes breaks into
two classes of homeomorphisms. The solution of the mystery was the incentive to produce this
paper. It reveals scientific collaboration at its best.

\section{Conclusion}
A closed orientable 3-manifold is denoted {\em $n$-small} if it is induced by surgery on
a blackboard framed link with at most $n$ crossings. We provide an instance of the general theory
to produce a recursive indexation of $n$-small 3-manifolds up to homeomorphism. We 
solve this problem
up to $n=9$. Conceptually we could go on forever, finding in the way 
tougher and tougher examples to be distinguished
by yet to be found new invariants. This is particularly reachable, if we consider that the 
census generating algorithms are fully parallellizable, and this parallelization was not used.
 The topological classification of the 9-small 3-manifolds involve three invariants: 
 $$INV=\{{\ homology}, \ W\hspace{-0.5mm}RT_{12}, { \ length \ of\ smallest\  geodesic}\ \}.$$
 The classification was nearly complete in \cite{lins2007blink}, except for two doubts. Recently,
 after we posted a challenge in the arXiv, \cite{linslins2013A} these doubts were solved by
 M. Culler, N. Dunfield and C. Hodgson using SnapPy \cite{snappy}. This made us add 
 $$ \ length \ of\ smallest\  geodesic $$ 
 which we define as 0, if the manifold is not hyperbolic, to our list of invariants. 
 The $9$-small 3-manifold classification maintains live the two Conjectures 
 of page 15 of \cite{lins1995gca} based on two kinds of moves $TS$ and $U$:
 the $TS$- and $U$-moves yield an efficient algorithm to classify 3-manifolds 
 by explicitly displaying homeomorphisms among them, whenever they exist. 
 A recent finding by C. Hodgson concerning manifolds $T[71]$ and $T[79]$ forming the 
 HG8QI-class $14_{24}^t$, in the notation of page 239 of \cite{lins2007blink} shows that 
 the 3 invariants are not enough to
 decide the pair. This pair is the first one of 11 pairs that we display as some tougher challenges
 to 3-manifold topologists, \cite{lins2013}. Hodgson's finding is that the volume as well as the
 lengths of the smallest geodesics fail to distinguish $T[71]$ and $T[79]$. He proves them 
 to be non-homeomorphic by more sofisticated techniques, involving drilling along the smallest geodesics
 to get non-isometric manifolds with toroidal boundary. Using SnapPy, GAP and Sage, 
 N. Dunfield shows that $T[71]$ and $T[79]$  are also distinguished by the homology of 
 its 5-covers. We had problems trying to find low index 
 subgroups for the $HG8QI$-class $15^t_{19}$. We did not put enough time for finding
 all the 6-covers. However, this now has been fixed.

\begin{figure}[H]
\begin{center}
\includegraphics[width=14cm] {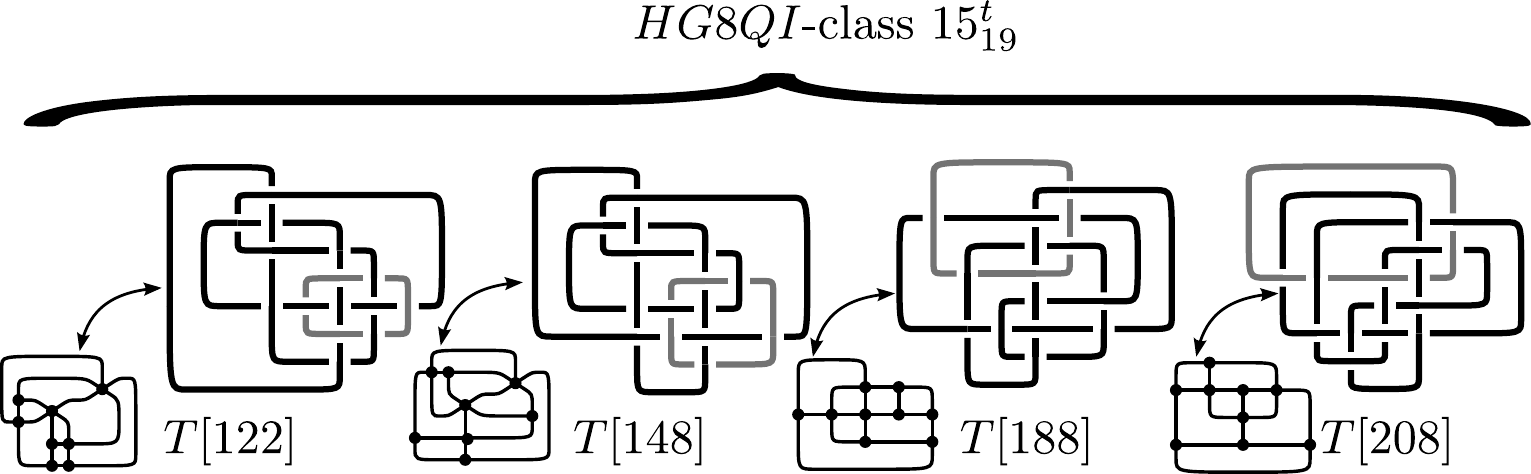}
\caption{\sf The third tougher challenge of \cite{lins2013}: 
the $HG8QI$-class $15^t_{19}$,
solved by M. Culler using SnapPy. 
It is easy to show that $T[122]=T[188]$ and that $T[148]=T[208]$ 
Enumeration of all the homomorphisms of the fundamental group onto $S_6$ shows that
$T[122]$ has 213 6-covers while $T[148]$ has ``only'' 211 6-covers. Again
the tantalizing fact that not only the $WRT_{12}$-invariants are the same but that 
also the volumes of 
$T[122]$ and $T[148]$ are the same. And, even so, the manifolds are non-homeomorphic!
A pattern has been created and we pose the question: do the WRT-invariants 
(as defined combinatorially in \cite{lins1995gca} and \cite{kauffman1994tlr}) define the
volume of a hyperbolic closed 3-manifold?
}
\label{fig:T15t_19focusedchallenge}
\end{center}
\end{figure}

\noindent
Here is how Marc Culler solved the $HG8QI$-class $15^t_{19}$, transcripting from
the blog on Lower Dimensional Topology:
\begin{verbatim}
Sóstenes, I dont know what you tried, but I did not have 
any particular trouble classifying the four manifolds in your 
family 15_19 last night. First I checked the volumes of the four links. 
The link complements for T122 and T188 both have volume 32.1904446? while those 
for T148 and T208 have volume 31.8897187? . 
Next I asked SnapPy to look for isometries between the pairs with the same 
volume. It quickly found framing-preserving homeomorphisms 
between the link complements in both pairs. 
So that shows that the closed manifold T122 is 
homeomorphic to T188, and T148 is homeomorphic to T208. 
Next I did the Dehn fillings, and noted that the volumes of the 
filled manifolds T122 and T148 were the same, namely 27.670218369?. 
To distinguish them, I asked SnapPy to find their 6-fold covers. 
This was taking a while, as one would expect with manifolds that large, 
so I went to bed. When I got up this morning SnapPy had found that 
T122 has 213 6-fold covers, one of which has first betti number 8. 
On the other hand, T148 has only 211 6-fold covers and the 
highest betti number that occurs is 7. That shows that the two 
closed manifolds are not homeomorphic. This family appears to 
behave very similarly to the others, except that the manifolds are larger. 
The same tools work fine to classify them, although it is true that 
the manifolds are approaching the size where computing the covering 
spaces becomes unpleasantly slow.
\end{verbatim} 

Using covering techniques, by now the majority of the 11 {\em Tough Challenges} of \cite{lins2013} has been solved.
So we have not yet arived to our currently technological limit.
The solution of these challenges depends on the existence of low index subgroups of the 
fundamental group of the manifolds. What about if these subgroups do not exist?
We feel that in the near future a software like 
BLINK can discover tougher and tougher
pairs of specific 3-manifolds to a point where itsef and no other software will be able to
prove equality or distinctiveness. 
At this point we will need to find algorithms for better invariants 
which perform quicker than the current ones. 
Invariants which are not so exponential
as looking for small index coverings which depend on the unpredictable 
idiosyncrasies of the fundamental group. The recent theorem proved that
closed oriented connected 3-manifolds are a subtle class of plane graphs (\cite{lins2013B}) 
may have something to say in regard to
these (yet to be found) invariants.

\bibliographystyle{plain}
\bibliography{bibtexIndex.bib}

\begin{thebibliography}{10}

\bibitem{ahuja1993network}
R.~K. Ahuja, T.~L. Magnanti, and J.~B. Orlin.
\newblock Network flows: theory, algorithms, and applications.
\newblock 1993.

\bibitem{alexandrov1961elementary}
P.~S. Aleksandroff.
\newblock {\em Elementary concepts of topology}, volume 747.
\newblock Courier Dover Publications, 1961.

\bibitem{alexandroff1932einfachste}
P.~Alexandroff and D.~Hilbert.
\newblock {\em Einfachste grundbegriffe der topologie}.
\newblock Springer Berlin, 1932.

\bibitem{bondy1976gta}
J.~A. Bondy and U.~S.~R. Murty.
\newblock {\em {Graph theory with applications}}.
\newblock Macmillan London, 1976.

\bibitem{snappy}
M.~Culler, N.~M. Dunfield, and J.~Weeks.
\newblock Snappy, a computer program for studying the geometry and topology of
  3-manifolds, 2012.

\bibitem{ferri1982crystallisation}
M.~Ferri and C.~Gagliardi.
\newblock {Crystallisation moves}.
\newblock {\em Pacific J. Math}, 100(1):85--103, 1982.

\bibitem{friedl3}
S.~Friedl.
\newblock 3-{M}anifolds after {P}erelman.

\bibitem{gap2002gap}
GAP Group et~al.
\newblock {GAP} -- {G}roups, {A}lgorithms, and {P}rogramming, version 4.3,
  2002, 2002.

\bibitem{hempel1976}
J.~Hempel.
\newblock {\em 3-Manifolds}, volume 349.
\newblock Amer, Math. Soc., 1976.

\bibitem{hodgson1994symmetries}
C.~D. Hodgson and J.~R. Weeks.
\newblock Symmetries, isometries and length spectra of closed hyperbolic
  three-manifolds.
\newblock {\em Experimental Mathematics}, 3(4):261--274, 1994.

\bibitem{kauffman1991knots}
L.H. Kauffman.
\newblock {\em Knots and physics}, volume~1.
\newblock World Scientific Publishing Company, 1991.

\bibitem{kauffman1994tlr}
L.H. Kauffman and S.~L. Lins.
\newblock {Temperley-Lieb Recoupling Theory and Invariants of 3-manifolds}.
\newblock {\em Annals of Mathematical Studies, Princeton University Press},
  134:1--296, 1994.

\bibitem{kirby1978calculus}
R.~Kirby.
\newblock A calculus for framed links in ${S}^3$.
\newblock {\em Inventiones Mathematicae}, 45(1):35--56, 1978.

\bibitem{Klarreich1012}
E.~Klarreich.
\newblock Getting into shapes: from hyperbolic geometry to cube complexes.
\newblock {\em Simons Foundation}, October, 2012.

\bibitem{lickorish1962representation}
W.~B.~R. Lickorish.
\newblock A representation of orientable combinatorial 3-manifolds.
\newblock {\em Annals of Mathematics}, 76(3):531--540, 1962.

\bibitem{lickorish1991three}
W.~B.~R. Lickorish.
\newblock Three-manifolds and the {T}emperley-{L}ieb algebra.
\newblock {\em Mathematische Annalen}, 290(1):657--670, 1991.

\bibitem{lins2007blink}
L.~D. Lins.
\newblock Blink: a language to view, recognize, classify and manipulate
  3{D}-spaces.
\newblock {\em Arxiv preprint math/0702057}, 2007.

\bibitem{lins2006blobs}
S.~Lins and M.~Mulazzani.
\newblock {Blobs and flips on gems}.
\newblock {\em Journal of Knot Theory and its Ramifications}, 15(8):1001--1035,
  2006.

\bibitem{lins1995gca}
S.~L. Lins.
\newblock {\em {Gems, Computers, and Attractors for 3-Manifolds}}.
\newblock World Scientific, 1995.

\bibitem{lins2013B}
S.~L. Lins.
\newblock Closed oriented 3-manifolds are equivalence classes of plane graphs.
\newblock {\em arXiv:1305.4540v3 [math.GT]}, 2013.

\bibitem{lins2013}
S.~L. Lins.
\newblock A tougher challenge to 3-manifold topologists and group algebraists.
\newblock {\em arXiv:1305.2617v2 [math.GT]}, 2013.

\bibitem{linslins2013A}
S.~L. Lins and L.~D. Lins.
\newblock A challenge to 3-manifold topologists and group algebraists.
\newblock {\em arXiv:1213.5964v4 [math.GT]}, 2013.

\bibitem{linsmachadoA2012}
S.~L. Lins and R.~N. Machado.
\newblock Framed link presentations of 3-manifolds by an ${O}(n^2)$ algorithm,
  {I}: gems and their duals.
\newblock {\em arXiv:1211.1953v2 [math.GT]}, 2012.

\bibitem{linsmachadoB2012}
S.~L. Lins and R.~N. Machado.
\newblock Framed link presentations of 3-manifolds by an ${O}(n^2)$ algorithm,
  {II}: colored complexes and boundings in their complexity.
\newblock {\em arXiv:1212.0826v2 [math.GT]}, 2012.

\bibitem{linsmachadoC2012}
S.~L. Lins and R.~N. Machado.
\newblock Framed link presentations of 3-manifolds by an ${O}(n^2)$ algorithm,
  {III}: geometric complex $\mathcal{H}_n^\star$ embedded into $\mathbb{R}^3$.
\newblock {\em arXiv:1212.0827v2 [math.GT]}, 2012.

\bibitem{lovasz1998om}
L.~Lovasz.
\newblock {One Mathematics}.
\newblock {\em The Berliner Intelligencer, Berlin}, pages 10--15, 1998.

\bibitem{magnus1966combinatorial}
W.~Magnus, A.~Karrass, and D.~Solitar.
\newblock {\em {Combinatorial group theory}}.
\newblock Interscience Publishers New York, 1966.

\bibitem{martelli2011finite}
B.~Martelli.
\newblock A finite set of local moves for {K}irby calculus.
\newblock {\em Arxiv preprint arXiv:1102.1288}, 2011.

\bibitem{martelli2012finite}
B.~Martelli.
\newblock A finite set of local moves for {K}irby calculus.
\newblock {\em Journal of Knot Theory and Its Ramifications}, 21(14), 2012.

\bibitem{matveev2007algorithmic}
S.~Matveev.
\newblock {\em Algorithmic topology and classification of 3-manifolds},
  volume~9.
\newblock Springer, 2007.

\bibitem{moise1952affine}
E.~E. Moise.
\newblock Affine structures in 3-manifolds: V. the triangulation theorem and
  hauptvermutung.
\newblock {\em The Annals of Mathematics}, 56(1):96--114, 1952.

\bibitem{reshetikhin1991invariants}
N.~Reshetikhin and V.~G. Turaev.
\newblock Invariants of 3-manifolds via link polynomials and quantum groups.
\newblock {\em Inventiones mathematicae}, 103(1):547--597, 1991.

\bibitem{rolfsen2003knots}
D.~Rolfsen.
\newblock {\em Knots and links}.
\newblock American Mathematical Society, 2003.

\bibitem{sage2012}
W.~Stein.
\newblock {\em {Sage}: {O}pen {S}ource {M}athematical {S}oftware ({V}ersion
  2.10.2)}.
\newblock The Sage~Group, 2008.
\newblock {\tt http://www.sagemath.org}.

\bibitem{tamassia1987egg}
R.~Tamassia.
\newblock {On Embedding a Graph in the Grid with the Minimum Number of Bends}.
\newblock {\em SIAM Journal on Computing}, 16:421, 1987.

\bibitem{turaev1992state}
V.~G. Turaev and O.~Y. Viro.
\newblock State sum invariants of 3-manifolds and quantum 6j-symbols.
\newblock {\em Topology}, 31(4):865--902, 1992.

\bibitem{weeks2001snappea}
J.~Weeks.
\newblock Snap{P}ea: a computer program for creating and studying hyperbolic
  3-manifolds, 2001.

\bibitem{witten1989quantum}
E.~Witten.
\newblock Quantum field theory and the {J}ones polynomial.
\newblock {\em Communications in Mathematical Physics}, 121(3):351--399, 1989.

\end{thebibliography}

\vspace{5mm}
\begin{center}
\hspace{7mm}
\begin{tabular}{l}
   S\'ostenes L. Lins\\
   Centro de Inform\'atica, UFPE \\
   Av. Jornalista Anibal Fernandes s/n\\
   Recife, PE 50740-560 \\
   Brazil\\
   sostenes@cin.ufpe.br
\end{tabular}
\hspace{20mm}
\hspace{7mm}
\begin{tabular}{l}
Lauro D. Lins\\
AT\&T Labs Research \\
180 Park Avenue \\
Florham Park, NJ 07932 \\
USA\\
llins@research.att.com
\end{tabular}
\hspace{20mm}

\end{center}
\newpage
\section{Appendix: census (no misses, no duplicates) of 9-small 3-manifolds}
\subsection*{Part 1/4 in terms of blinks:}
\center{\includegraphics[width=16.5cm]{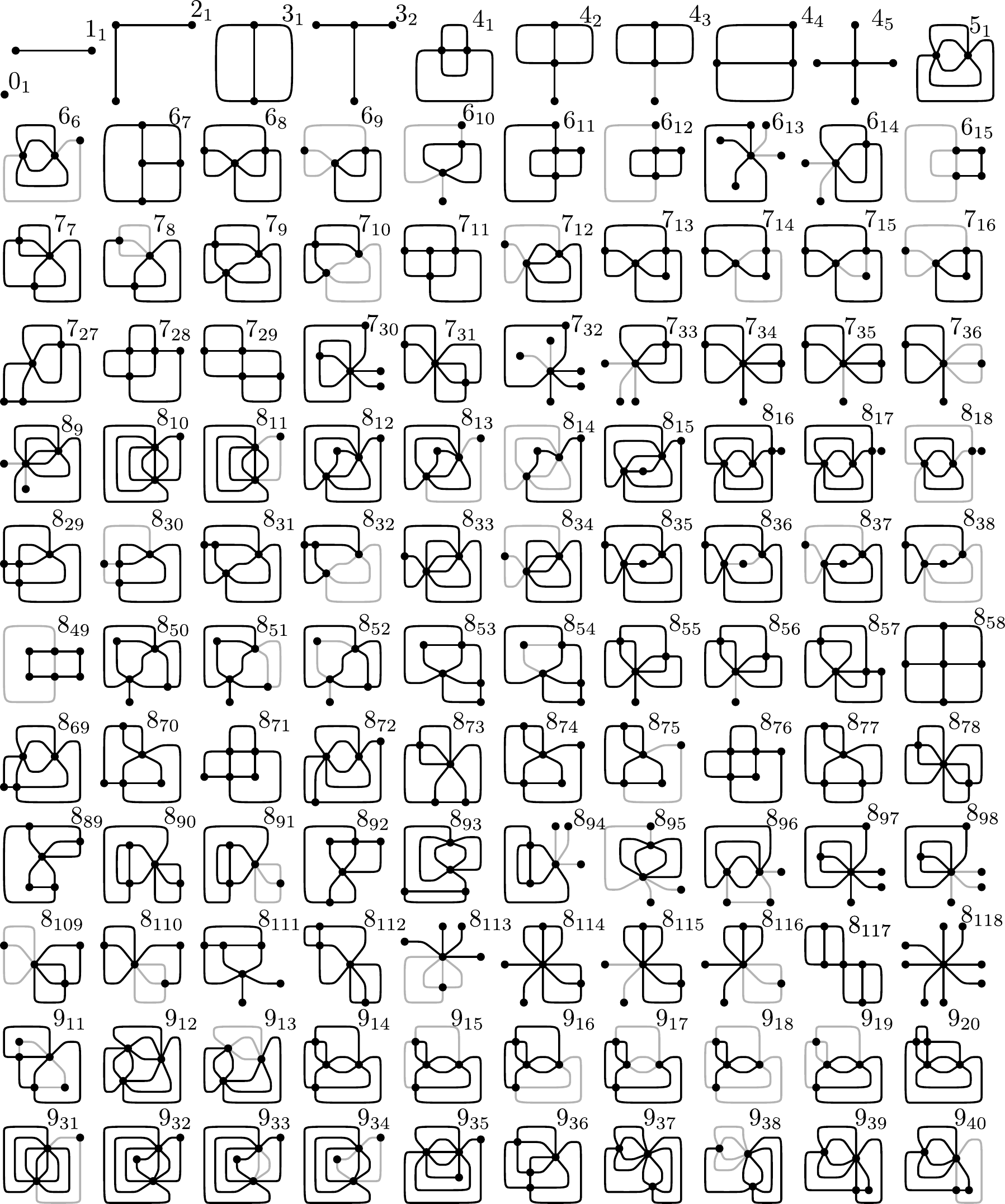}} \eject \noindent
\subsection*{Part 1/4 in terms of blackboard framed links:}
\center{\includegraphics[width=16.5cm]{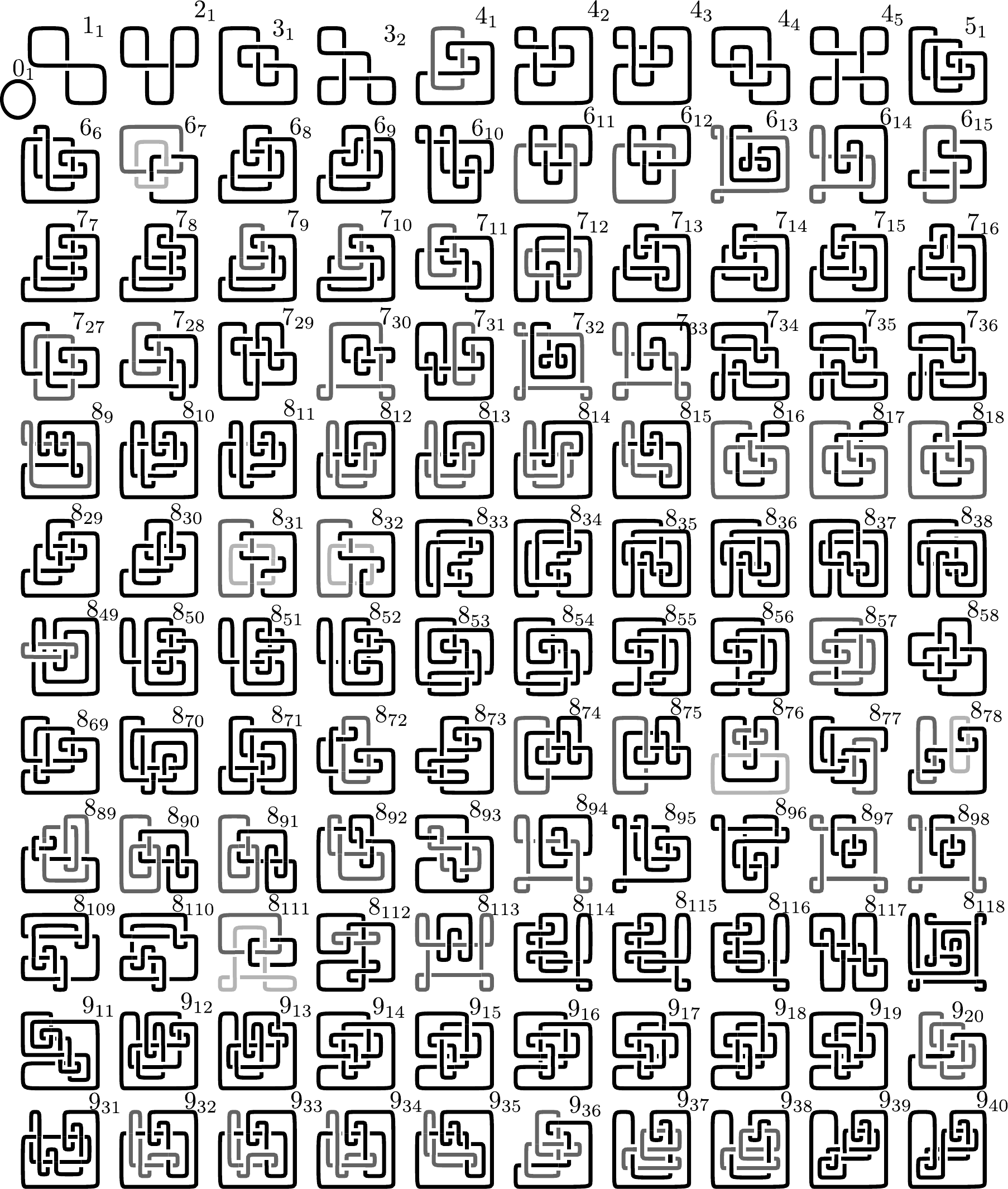}} \eject
\subsection*{Part 2/4 in terms of blinks:}
\center{\includegraphics[width=16.5cm]{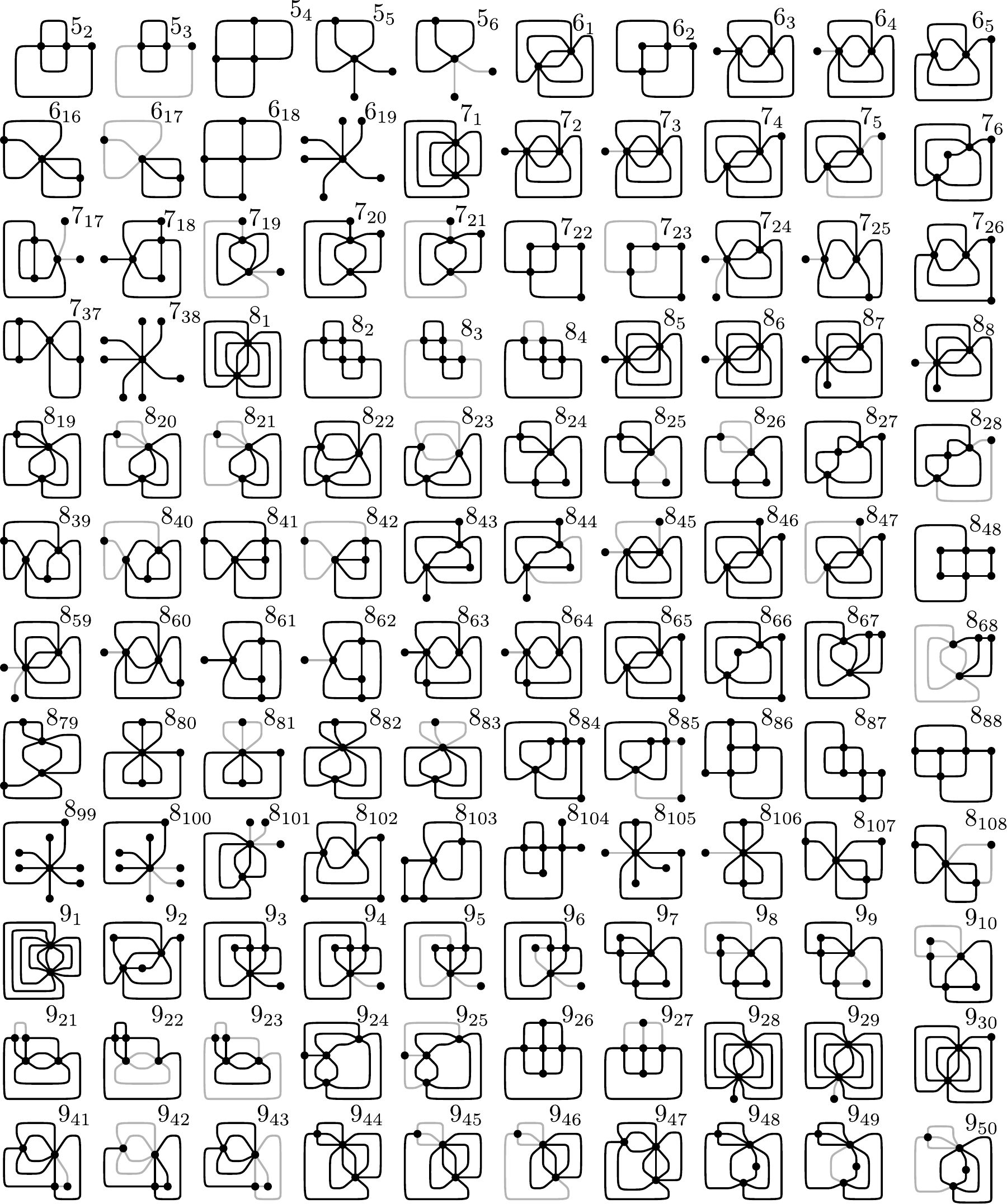}} \eject
\subsection*{Part 2/4 in terms of blackboard framed links:}
\center{\includegraphics[width=16.5cm]{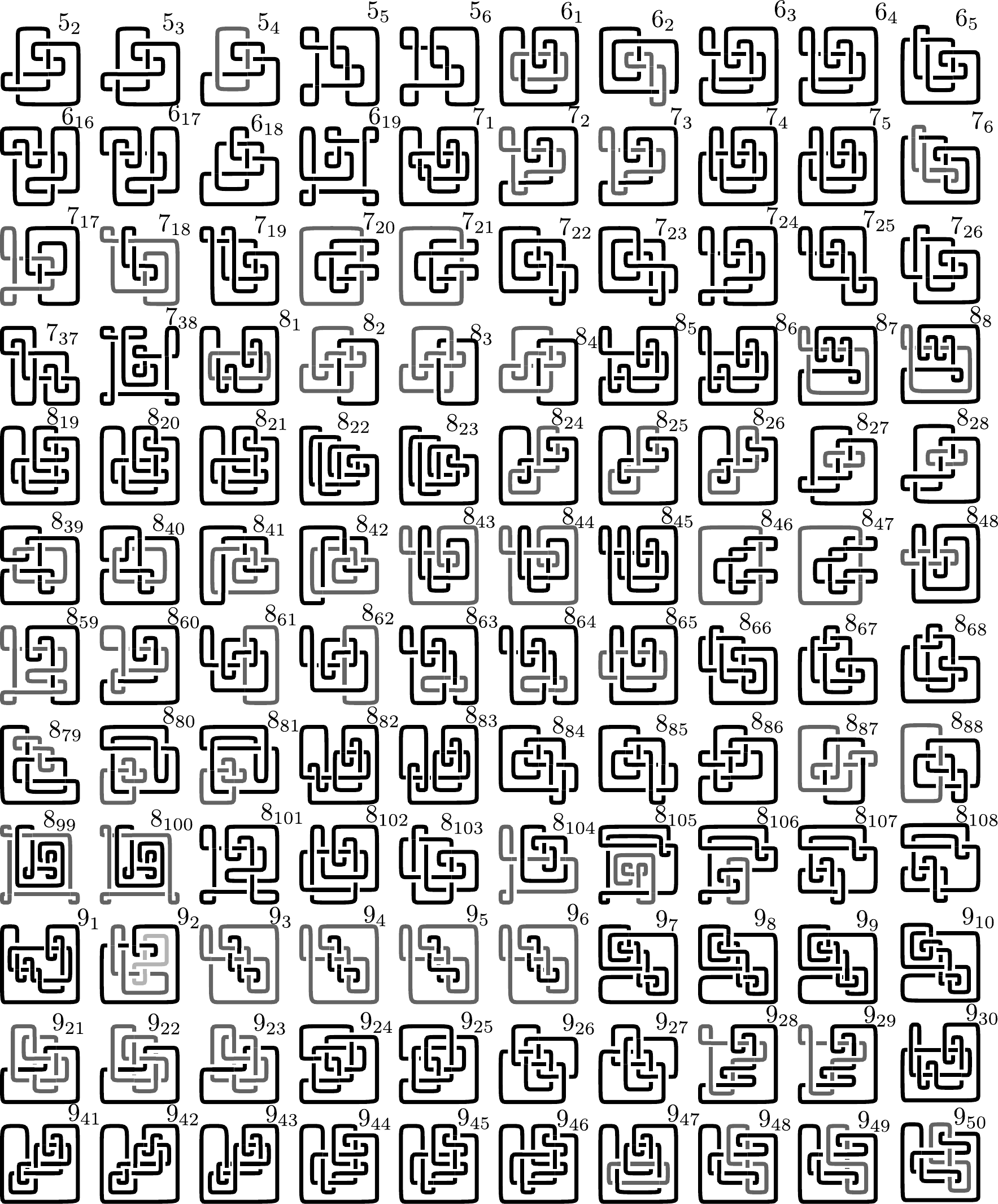}} \eject
\subsection*{Part 3/4 in terms of blinks:}
\center{\includegraphics[width=16.5cm]{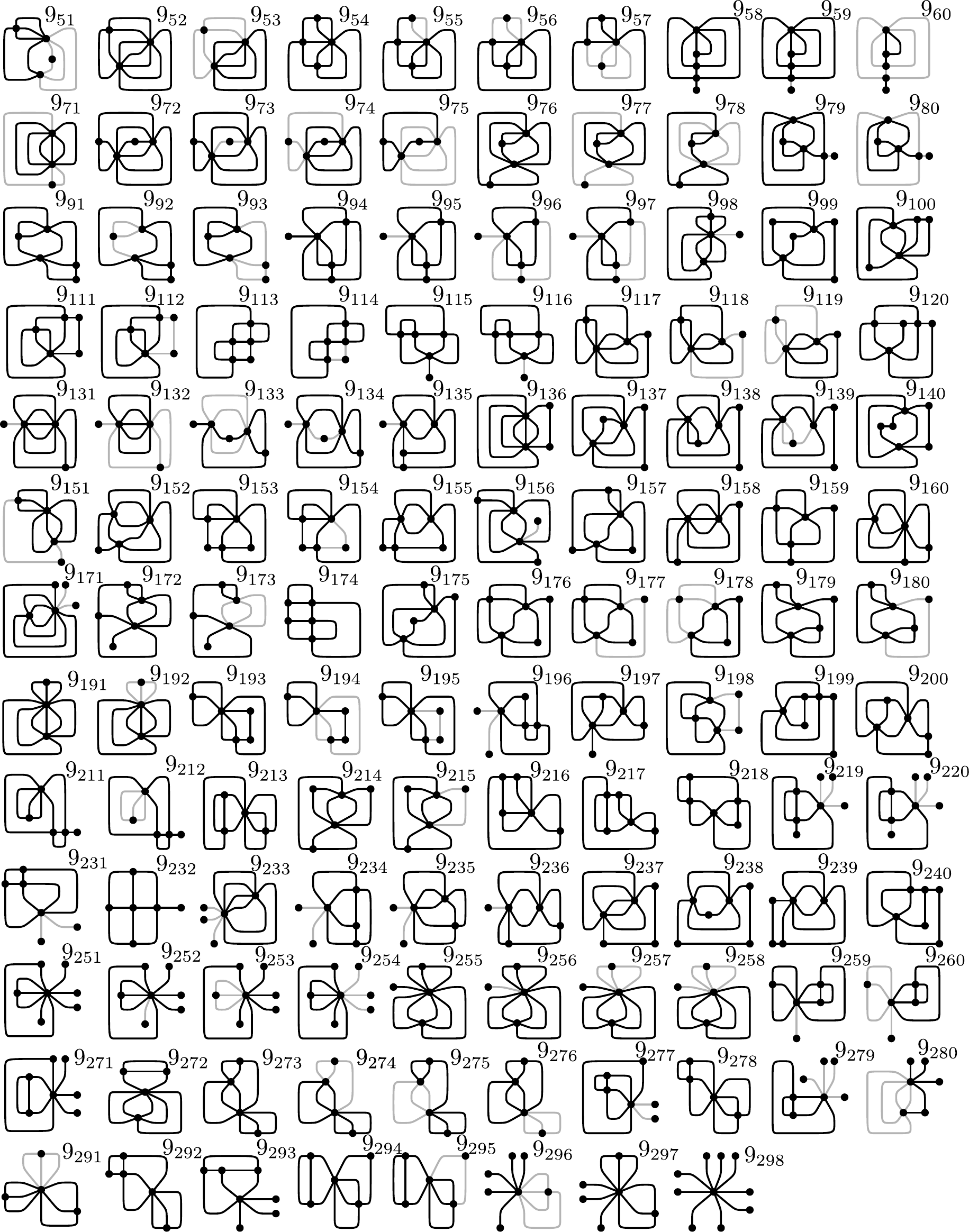}} \eject
\subsection*{Part 3/4 in terms of blackboard framed links:}
\center{\includegraphics[width=16.5cm]{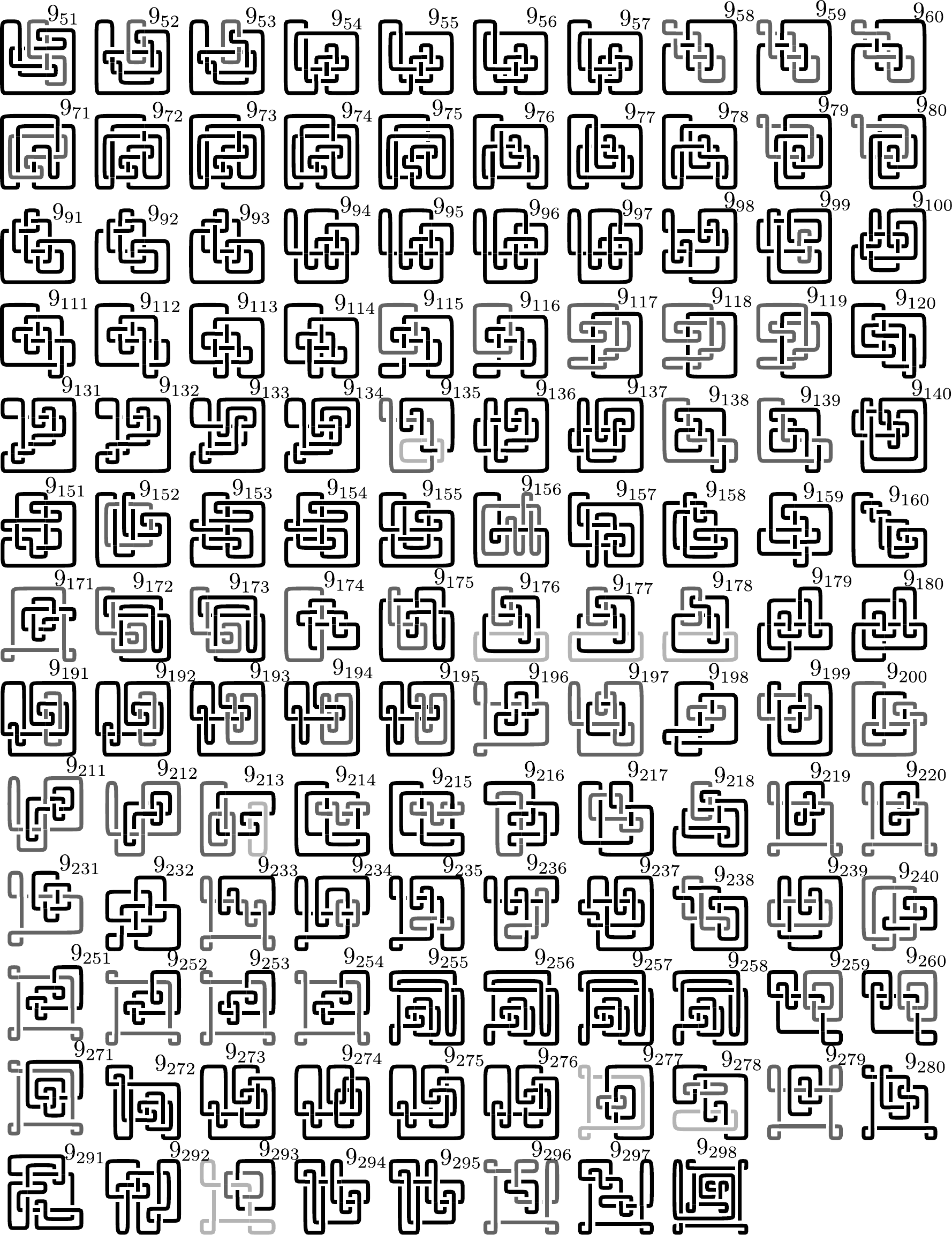}} \eject
\subsection*{Part 4/4 in terms of blinks:}
\center{\includegraphics[width=16.5cm]{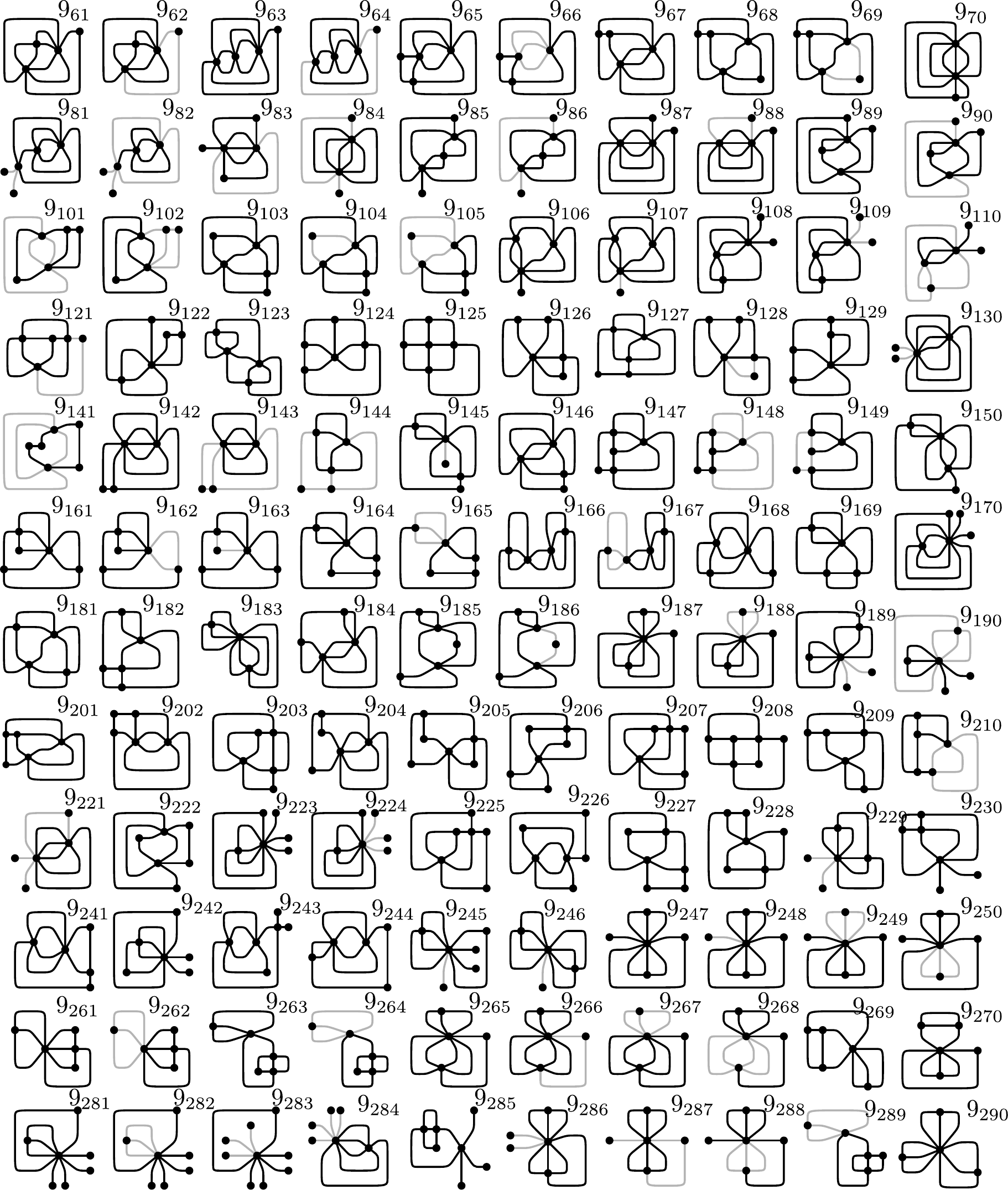}} \eject
\subsection*{Part 4/4 in terms of blackboard framed links:}
\center{\includegraphics[width=16.5cm]{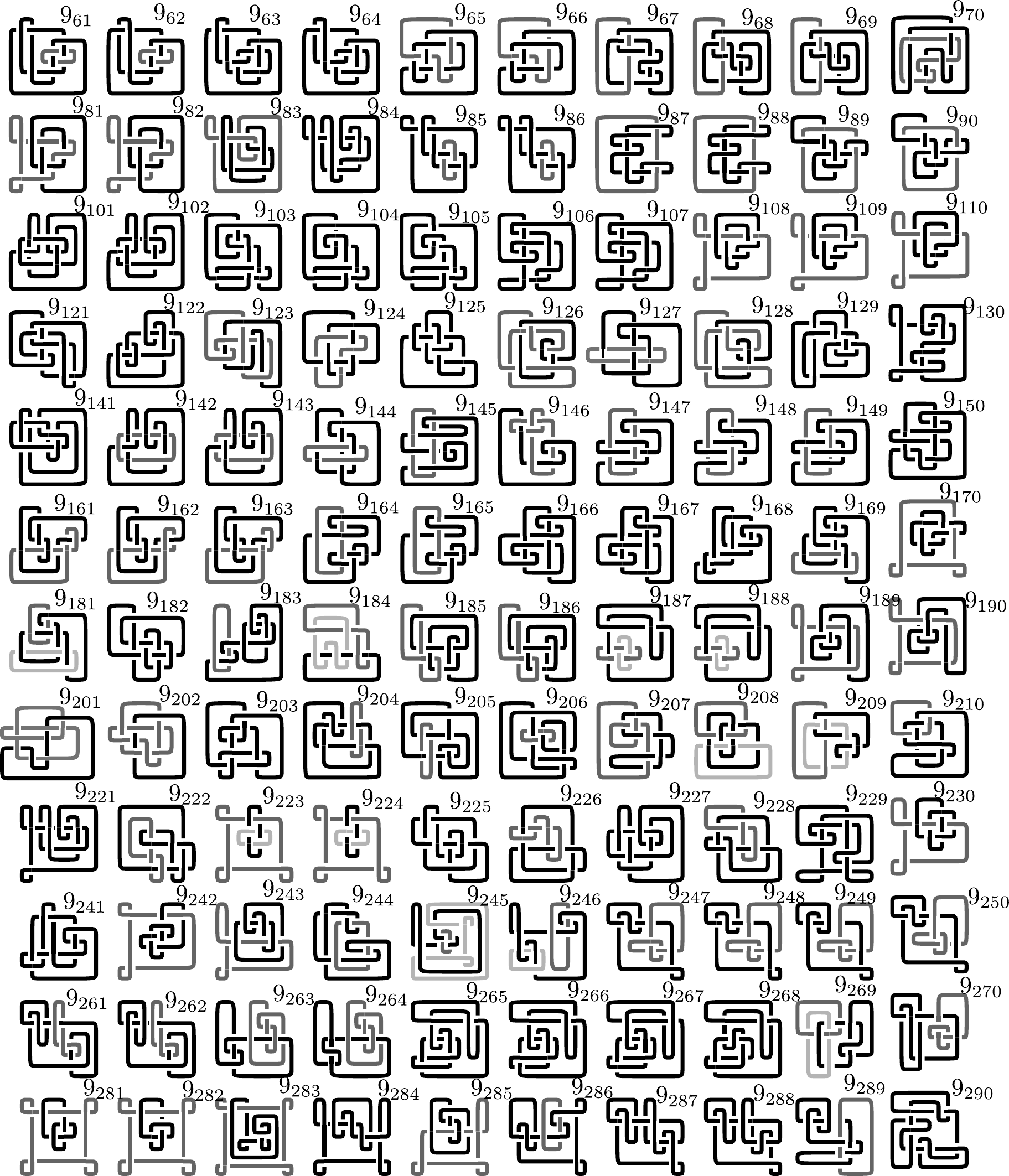}}\eject

\end{document}